\documentclass{mcom-l}

\usepackage{amsmath,amssymb,amsfonts,mathtools}
\usepackage{bm}
\usepackage{booktabs}
\usepackage{array}
\usepackage{enumitem}
\usepackage{microtype}

\usepackage{url}

\newtheorem{theorem}{Theorem}[section]
\newtheorem{proposition}[theorem]{Proposition}
\newtheorem{lemma}[theorem]{Lemma}
\newtheorem{corollary}[theorem]{Corollary}
\newtheorem{remark}[theorem]{Remark}
\newtheorem{assumption}[theorem]{Assumption}
\theoremstyle{definition}
\newtheorem{definition}[theorem]{Definition}

\newcommand{\R}{\mathbb{R}}

\newcommand{\D}{\mathcal{D}}
\newcommand{\calR}{\mathcal{R}}
\newcommand{\calS}{\mathcal{S}}

\newcommand{\inner}[2]{(#1,#2)}
\newcommand{\norm}[1]{\lVert #1\rVert}
\newcommand{\abs}[1]{\lvert #1\rvert}

\newcommand{\bdfthree}{D_3}
\newcommand{\bdftwo}{D_2}
\newcommand{\be}{D_1}

\newcommand{\Lamthree}{\Lambda_3}
\newcommand{\LamthreeMu}{\Lambda_3(\mu)}
\newcommand{\betabar}{\overline\beta}

\title[Long-time BDF3 method for parabolic equation with fading memory]
{Long-time stability and convergence of a BDF3 method for a parabolic equation with fading memory}

\author{Wenlin Qiu}
\address{School of Mathematics, Yunnan Normal University, Kunming 650500, P. R. China.}
\email{qwllkx12379@163.com}

\subjclass[1991]{Primary 65M12, 65M15, 65R20; Secondary 45K05, 65N30}
\keywords{parabolic integro-differential equation, fading memory, BDF3, long-time stability, uniform error estimate, finite element}

\begin{document}

\begin{abstract}
We analyze the long-time stability and convergence of a third-order backward differentiation formula (BDF3) for parabolic equations with exponentially fading memory.  The BDF3 scalar resolvent has oscillatory coefficients, so the monotonicity argument used for lower-order formulas in [Math. Comp. 62 (1994) 477-496] does not apply. We instead control memory feedback by the absolute mass of its discrete Green sequence and establish the rigorous stability-constant bounds $23/20<\Lambda_3<579/500$.  A causal third-order memory quadrature with nonnegative assembled weights then yields exponentially weighted stability without a parabolic step-size restriction.  With third-order starting values, we prove a uniform-in-time third-order temporal error estimate and extend it to conforming finite elements. Numerical tests confirm the predicted convergence orders and the long-time behavior.
\end{abstract}

\maketitle

\section{Introduction}

We consider the linear parabolic integro-differential problem
\begin{equation}\label{eq:model}
 u_t(t)+Au(t)=\int_0^t b(t-s)Au(s)\,ds+f(t),\qquad t>0,\qquad u(0)=u_0,
\end{equation}
posed in a real Hilbert space $H$.  The operator $A$ is assumed to be selfadjoint and strictly positive with compact inverse, while the scalar kernel $b$ decays exponentially.  Problems of the form \eqref{eq:model} arise, for example, in heat conduction with memory and in other evolution models with hereditary effects; see \cite{ColemanGurtin1967,SloanThomee1986,ThomeeWahlbin1994} and the references therein.

The issue addressed in this paper is not merely convergence on a prescribed finite interval.  Our aim is to derive estimates that remain informative when the final time becomes very large.  A conventional discrete Gronwall argument typically introduces a factor that grows exponentially with the final time and therefore obscures the dissipative character of the continuous equation.  For fading-memory problems this is particularly undesirable: the continuous solution may decay exponentially while a finite-time error estimate becomes increasingly pessimistic as the computation proceeds.

A systematic long-time analysis for backward Euler and the second-order backward differentiation formula was developed by Thom\'ee-Wahlbin \cite{ThomeeWahlbin1994}.  In the special but important case in which the same elliptic operator appears in the instantaneous and memory terms, their argument diagonalizes the operator and reduces the stability question to scalar recurrences.  For backward Euler the scalar resolvent is monotone, whereas BDF2 introduces two characteristic roots and a stability constant measuring the absolute mass of the corresponding discrete Green function.  This makes the BDF2 admissible memory range slightly smaller than its continuous counterpart.  Finite-time Galerkin analyses of related parabolic Volterra equations, including the use of Ritz--Volterra projections, may be found in \cite{LinThomeeWahlbin1991,PaniThomeeWahlbin1992,ThomeeZhang1989,Zhang1993}.  The convolution-quadrature viewpoint developed by Lubich \cite{Lubich1988I,Lubich1988II} provides a different and very influential route to discretizing temporal convolutions.

High-order BDF stability analysis has advanced through energy estimates, multipliers, and discrete gradient structures; see \cite{AkrivisChenYu2025,AkrivisChenYuZhou2021,ChenShiYuZhou2026,ContriKovacsMassing2025,HuangShen2024,LiaoKang2024,WuZhou2021}.  Related high-order methods for nonlocal evolution equations include exponential or BDF convolution quadrature, rational kernel approximation, and spectral treatments of delay terms \cite{JinQuanWohlmuthZhou2025,KhristenkoWohlmuth2023,LiLinMaRao2023,LiMa2022,LiuMa2025,LopezFernandezPlociniczak2025,ShiChen2023,ShiDuChen2024}.  Many of these analyses address singular or algebraically decaying kernels rather than exponentially fading memory.

Extending the long-time BDF2 results of \cite{ThomeeWahlbin1994} to BDF3 presents three obstacles.  Although the negative real axis lies in its stability region, BDF3 is not A-stable \cite{Dahlquist1963,HairerWanner1996}; its three characteristic roots produce an oscillatory Green sequence.  Moreover, third-order accuracy requires a matching memory quadrature and accurate causal starting values.

We address these difficulties by introducing the BDF3 scalar Green sequence associated with a spectral parameter $\mu=k\lambda$.  Using the Schur--Cohn criterion \cite{Cohn1922,Schur1918}, we prove that all three roots lie strictly inside the unit disk for every $\mu>0$.  We then define
\begin{equation}\label{eq:Lambda3intro}
 \Lamthree:=\sup_{\mu>0}\frac{6\mu}{11+6\mu}\sum_{m=0}^{\infty}\abs{G_m(\mu)},
\end{equation}
where $G_m(\mu)$ is the BDF3 discrete Green function.  This extends the BDF2 stability constant of \cite{ThomeeWahlbin1994}.  We prove that it is finite and satisfies
\[
 \frac{23}{20}<\Lamthree<\frac{579}{500}.
\]
We obtain the upper bound by root parametrization and exact Bernstein enclosures with controlled tails.  Numerically, $\Lamthree\approx1.1501238$ at $\mu\approx7.43$.  The analytic enclosure yields the sufficient kernel condition $\beta<500/579$.

We combine a Simpson starting block with Adams--Bashforth quadrature.  The assembled weights are nonnegative, which supports the weighted-memory argument.  First- and third-derivative quadrature estimates provide stability and consistency, respectively.

The analysis has three parts.  First, we prove Green-function bounds uniformly over the spectrum of $A$, and bound the amplification factor $\Lambda_3$.  Next, we combine these bounds with the nonnegative third-order quadrature weights to obtain exponentially weighted stability and a uniform $O(k^3)$ error estimate.  Finally, we derive a mesh-independent estimate for conforming finite elements.  The appendices contain the exact enclosure of $\Lambda_3$ and the technical weighted-resolvent estimates.

\section{Problem setting}\label{sec:continuous}

Let $H$ be a real Hilbert space with inner product $\inner{\cdot}{\cdot}$ and norm $\norm{\cdot}$.  We assume that $A:\D(A)\subset H\to H$ is selfadjoint, strictly positive, and has compact inverse.  Hence there is an orthonormal basis $\{\varphi_m\}_{m\ge1}$ and a nondecreasing sequence of positive eigenvalues such that
\[
 A\varphi_m=\lambda_m\varphi_m,\qquad 0<\lambda_0\le \lambda_1\le\lambda_2\le\cdots,\qquad \lambda_m\to\infty.
\]
For $r\in\R$ we use the scale norm $\norm{v}_r=\norm{A^{r/2}v}$ whenever the right-hand side is defined.

We take $b(t)=0$ for $t<0$.  The kernel may change sign; our stability estimates therefore use its absolute mass.

We use the following decay condition, imposing more derivatives only for convergence.

\begin{definition}\label{def:kernel}
Let $q$ be a nonnegative integer.  We say that the causal kernel $b$ has exponentially fading memory of order $q$ if $b\in C^q([0,\infty))$ and there are constants $\alpha>0$ and $K_j\ge0$, $0\le j\le q$, for which
\[
 \abs{b^{(j)}(t)}\le K_j e^{-\alpha t},\qquad t\ge0,\qquad 0\le j\le q.
\]
Here and below $b^{(0)}=b$.  We denote the absolute memory mass by
\begin{equation}\label{eq:beta}
 \beta:=\int_0^\infty\abs{b(s)}\,ds,
\end{equation}
and, for $0\le\gamma<\alpha$, the exponentially weighted mass by
\begin{equation}\label{eq:beta-gamma}
 \beta_\gamma:=\int_0^\infty e^{\gamma s}\abs{b(s)}\,ds.
\end{equation}
\end{definition}

\begin{assumption}\label{ass:kernel}
For the long-time stability analysis we assume that $b$ has exponentially fading memory of order one.  For the third-order convergence analysis we assume order three.  In both cases the absolute mass \eqref{eq:beta} is finite.
\end{assumption}

The BDF3 stability theory will use the structural smallness condition
\begin{equation}\label{eq:kernel-sharp-condition}
 \Lamthree\beta<1.
\end{equation}
The constant $\Lamthree$ is defined in Section~\ref{sec:resolvent}.  Theorem~\ref{thm:Lambda-enclosure} will show that the directly checkable condition
\begin{equation}\label{eq:kernel-explicit-condition}
 \beta<\frac{500}{579}
\end{equation}
is sufficient for \eqref{eq:kernel-sharp-condition}.

Decay of $b$ and $b'$ controls the weighted discrete memory mass.  The third-order consistency estimate additionally uses $b''$ and $b^{(3)}$.

The class above includes the standard kernels used in heat conduction with memory.  For example, if
\[
 b(t)=\kappa e^{-at},\qquad a>0,
\]
then $\beta=\abs{\kappa}/a$, and the explicit BDF3 condition \eqref{eq:kernel-explicit-condition} reduces to $\abs{\kappa}/a<500/579$.  More generally, for a finite Prony series
\[
 b(t)=\sum_{\ell=1}^L\kappa_\ell e^{-a_\ell t},\qquad a_\ell>0,
\]
we have $\beta\le\sum_{\ell=1}^L\abs{\kappa_\ell}/a_\ell$, so the latter sum being smaller than $500/579$ is sufficient.  Damped oscillatory kernels, such as $b(t)=\kappa e^{-at}\cos(\omega t)$, are also covered because $\beta\le\abs{\kappa}/a$.  These examples emphasize that our theory is based on fading amplitude rather than positivity of the kernel.

The continuous equation already reveals the central role of $\beta$.  If $f=0$ and $v(t)=e^{\gamma t}u(t)$, then $v$ satisfies the same type of equation with the modified kernel $e^{\gamma t}b(t)$ and with $A$ replaced by $A-\gamma I$.  Choosing $\gamma>0$ sufficiently small and using $\lambda_0\norm{v}^2\le\inner{Av}{v}$ gives the standard fading-memory estimate.  We record it for reference.

\begin{proposition}\label{prop:continuous}
Let $b$ have exponentially fading memory of order zero and assume $\beta<1$.  Then there are $\gamma>0$ and $C>0$ such that the solution of \eqref{eq:model} satisfies
\begin{equation}\label{eq:continuous-stability}
 \norm{u(t)}\le C e^{-\gamma t}\norm{u_0}+C\int_0^t e^{-\gamma(t-s)}\norm{f(s)}\,ds,
\end{equation}
for all $t\ge0$.
\end{proposition}

The estimate \eqref{eq:continuous-stability} is classical; see, for example, \cite{ThomeeWahlbin1994}.  Its constant is independent of the final time.  Since $\Lamthree\ge1$, the BDF3 condition \eqref{eq:kernel-sharp-condition} automatically implies $\beta<1$ and is therefore consistent with the continuous stability threshold.

\section{BDF3 discretization}\label{sec:scheme}

Let $k>0$ be a uniform time step and $t_n=nk$.  For a sequence $\{V^n\}$ we define the first three backward differentiation operators by
\[
 \be V^n=\frac{V^n-V^{n-1}}{k},\qquad
 \bdftwo V^n=\frac{3V^n-4V^{n-1}+V^{n-2}}{2k},
\]
and
\begin{equation}\label{eq:BDF3}
 \bdfthree V^n=\frac{11V^n-18V^{n-1}+9V^{n-2}-2V^{n-3}}{6k}.
\end{equation}
The main time-stepping method will be used for $n\ge3$.

\subsection{The memory quadrature}

For a scalar- or $H$-valued function $g$ sampled at $t_j$, we define the causal third-order quadrature for $n\ge3$ by
\begin{equation}\label{eq:quadrature-local}
 \sigma^n(g)=\frac{k}{3}\bigl(g^0+4g^1+g^2\bigr)
 +k\sum_{\ell=2}^{n-1}\left(\frac{23}{12}g^\ell-\frac{4}{3}g^{\ell-1}+\frac{5}{12}g^{\ell-2}\right).
\end{equation}
The first term is Simpson's rule on $[0,2k]$, while the subsequent intervals are treated by the classical three-step Adams--Bashforth formula \cite{HairerNorsettWanner1993}. Since the last value used in \eqref{eq:quadrature-local} is $g^{n-1}$, the resulting quadrature is causal.

Although the local Adams--Bashforth formula contains a negative coefficient, cancellation between consecutive intervals produces nonnegative global weights.  In particular, for $n\ge5$, \eqref{eq:quadrature-local} is equivalent to
\begin{equation}\label{eq:quadrature-global}
 \sigma^n(g)=k\left(\frac34g^0+\frac{5}{12}g^1+\frac43g^2+\sum_{j=3}^{n-3}g^j+\frac{7}{12}g^{n-2}+\frac{23}{12}g^{n-1}\right),
\end{equation}
where an empty sum is omitted.  For the two exceptional cases one obtains
\[
 \sigma^3(g)=k\left(\frac34g^0+\frac94g^2\right),\qquad
 \sigma^4(g)=k\left(\frac34g^0+\frac{5}{12}g^1+\frac{11}{12}g^2+\frac{23}{12}g^3\right).
\]
Thus every assembled quadrature weight is nonnegative. We denote the quadrature error by
\begin{equation}\label{eq:qdef}
 q^n(g):=\sigma^n(g)-\int_0^{t_n}g(s)\,ds.
\end{equation}
For the error defined in \eqref{eq:qdef}, the following estimate will be used twice, with different derivative orders.

\begin{lemma}\label{lem:quadrature}
Let $n\ge3$.  If $g\in W^{j,1}(0,t_n;H)$ for $j=1$ or $j=3$, then
\begin{equation}\label{eq:quadrature-error}
 \norm{q^n(g)}\le Ck^j\int_0^{t_n}\norm{g^{(j)}(s)}\,ds,
\end{equation}
where $C$ is independent of $k$, $n$, and $g$.
\end{lemma}
\begin{proof}
  A proof is given in Appendix~\ref{app:quadrature}.  The case $j=1$ is needed only to compare the discrete weighted mass of the kernel with the integral in \eqref{eq:beta}; the case $j=3$ gives the third-order consistency of the memory approximation.
\end{proof}

For the convolution in \eqref{eq:model}, we write $g_n(s)=b(t_n-s)Au(s)$ for $0\le s\le t_n$, and use $\sigma^n(g_n)$ in the scheme.

\subsection{The BDF3 scheme and starting values}

With the quadrature \eqref{eq:quadrature-local}, our main method is
\begin{equation}\label{eq:scheme}
 \bdfthree U^n+AU^n=\sigma^n(b_nAU)+f^n,\qquad n\ge3,
\end{equation}
where $f^n=f(t_n)$ and
\[
 \sigma^n(b_nAU)=\sigma^n\bigl(\{b(t_n-t_j)AU^j\}_{j=0}^{n-1}\bigr).
\]
Because the quadrature does not use $U^n$, each step requires one solve with the operator $11(6k)^{-1}I+A$.

A three-step formula requires two starting values.  We use a Crank--Nicolson step followed by one BDF2 step:
\begin{equation}\label{eq:start1}
 \be U^1+\frac12A(U^1+U^0)=\frac{k}{2}b(k)AU^0+\frac12(f^1+f^0),\qquad U^0=u_0,
\end{equation}
and
\begin{equation}\label{eq:start2}
 \bdftwo U^2+AU^2=2k\,b(k)AU^1+f^2.
\end{equation}
The memory contribution in \eqref{eq:start1} is obtained by averaging the equations at $t_0$ and $t_1$ and approximating the integral at $t_1$ by a left rectangle.  The memory contribution in \eqref{eq:start2} is the midpoint approximation to the integral over $[0,2k]$.  Both steps are causal.  Their accuracy is analyzed in Section~\ref{sec:convergence}.

\begin{remark}
The stability proof applies to any causal quadrature with nonnegative weights and sufficiently small weighted memory mass.  Third-order convergence additionally requires \eqref{eq:quadrature-error} with $j=3$ and $O(k^3)$ starting errors.
\end{remark}

\section{Scalar BDF3 resolvent}\label{sec:resolvent}

We analyze the BDF3 recurrence on each eigenmode of $A$.  Let $\lambda>0$ be an eigenvalue of $A$ and set
\[
 \mu=k\lambda,\qquad d_\mu=11+6\mu.
\]
After multiplying \eqref{eq:scheme} by $k$ and projecting onto the corresponding eigenvector, the homogeneous BDF3 operator takes the form
\begin{equation}\label{eq:scalar-rec}
 c^n-a_1(\mu)c^{n-1}+a_2(\mu)c^{n-2}-a_3(\mu)c^{n-3}=\eta(\mu)\,r^n,
\end{equation}
where
\begin{equation}\label{eq:coeffs}
 a_1(\mu)=\frac{18}{d_\mu},\qquad a_2(\mu)=\frac{9}{d_\mu},\qquad a_3(\mu)=\frac{2}{d_\mu},\qquad \eta(\mu)=\frac{6\mu}{d_\mu}.
\end{equation}
Here $r^n$ is the modal memory quadrature without the factor $\lambda$.  A projected forcing term adds $6k f_\lambda^n/d_\mu$ on the right-hand side.

\subsection{Root location}

The characteristic polynomial associated with the left-hand side of \eqref{eq:scalar-rec} is
\begin{equation}\label{eq:charpoly}
 p_\mu(r)=r^3-\frac{18}{d_\mu}r^2+\frac{9}{d_\mu}r-\frac{2}{d_\mu}.
\end{equation}
Although BDF3 is not A-stable, the negative real spectrum generated by $-A$ ensures the following root condition.

\begin{lemma}\label{lem:schur}
For every $\mu>0$, all three roots of $p_\mu$ satisfy $\abs{r}<1$.
\end{lemma}

\begin{proof}
We set $d=d_\mu>11$ and write
\[
 p(r)=r^3+ar^2+br+c,
 \qquad a=-\frac{18}{d},\quad b=\frac9d,\quad c=-\frac2d.
\]
Clearly $\abs{c}<1$.  For the reversed polynomial
\[
 p^*(r):=r^3p(r^{-1})=1+ar+br^2+cr^3,
\]
the first Schur transform is
\[
 \mathcal{S}p(r):=\frac{p(r)-cp^*(r)}{r}
 =(1-c^2)r^2+(a-bc)r+(b-ac).
\]
Substituting the above values of $a,b,c$ gives
\[
 \mathcal{S}p(r)=\frac1{d^2}\Big((d^2-4)r^2+18(1-d)r+9(d-4)\Big).
\]
Hence, after multiplication by the harmless positive factor $d^2$, the Schur transform is
\begin{equation}\label{eq:schurquadratic}
 Q_d(r)=(d^2-4)r^2+18(1-d)r+9(d-4).
\end{equation}
The Schur--Cohn reduction \cite{Cohn1922,Schur1918,Marden1966} states that, because $\abs{c}<1$, all zeros of $p$ lie in the open unit disk if and only if all zeros of $\mathcal{S}p$ do.  Thus it remains to verify that the two zeros of $Q_d$ lie in the unit disk.

For a real quadratic, let $Q(r)=A_{\diamond}r^2+B_{\diamond} r+C_{\diamond}$ with $A_{\diamond} >0$.  After normalization by $A_{\diamond}$, the second Schur reduction shows that its two zeros lie in the open unit disk precisely when
\[
 \abs{C_{\diamond}}<A_{\diamond},\qquad A_{\diamond}+B_{\diamond}+C_{\diamond}>0,\qquad A_{\diamond}-B_{\diamond}+C_{\diamond}>0.
\]
Indeed, if $q(r)=r^2+\alpha r+\beta$ with $\alpha=B_{\diamond}/A_{\diamond}$ and $\beta=C_{\diamond}/A_{\diamond}$, then $\abs{\beta}<1$ and
\[
 \frac{q(r)-\beta q^*(r)}{r}
 =(1-\beta)\big((1+\beta)r+\alpha\big).
\]
The remaining zero therefore lies in the unit disk exactly when
$\abs{\alpha}<1+\beta$, which is equivalent to the last two inequalities above.

For \eqref{eq:schurquadratic}, $A_{\diamond}=d^2-4>0$ and $C_{\diamond}=9(d-4)>0$.  Moreover,
\[
 A_{\diamond}-C_{\diamond}=d^2-4-9(d-4)=d^2-9d+32>0,
\]
where the quadratic on the right has negative discriminant, while
\[
 A_{\diamond}+B_{\diamond}+C_{\diamond}=(d^2-4)+18(1-d)+9(d-4)=(d-11)(d+2)>0
\]
and
\[
 A_{\diamond}-B_{\diamond}+C_{\diamond}=(d^2-4)-18(1-d)+9(d-4)=d^2+27d-58>0.
\]
All three conditions hold for $d>11$.  The Schur reduction therefore yields the desired conclusion.
\end{proof}

\subsection{The discrete Green function}

The discrete Green sequence is the unit-impulse response of the
homogeneous scalar difference operator, with zero prehistory.

We define the BDF3 Green sequence by
\begin{equation}\label{eq:Greenrec}
 G_0(\mu)=1,\qquad G_1(\mu)=a_1(\mu),\qquad G_2(\mu)=a_1(\mu)^2-a_2(\mu),
\end{equation}
with
\[
 G_m(\mu)=a_1(\mu)G_{m-1}(\mu)-a_2(\mu)G_{m-2}(\mu)+a_3(\mu)G_{m-3}(\mu),\qquad m\ge3.
\]
The associated generating function and signed-mass identity are
\begin{equation}\label{eq:GreenGF}
 \sum_{m=0}^{\infty}G_m(\mu)z^m=\frac{1}{1-a_1(\mu)z+a_2(\mu)z^2-a_3(\mu)z^3}
\end{equation}
and
\begin{equation}\label{eq:signedmass}
 \eta(\mu)\sum_{m=0}^{\infty}G_m(\mu)=1.
\end{equation}
The derivation, together with the relation between the signed and absolute Green masses, is recorded below.

For zero prehistory, the recurrence $L_\mu c^n=\eta(\mu)r^n$ has the convolution representation
\[
 c^n=\eta(\mu)\sum_{j=0}^nG_{n-j}(\mu)r^j.
\]
Thus $\Lambda_3(\mu)$ measures its $\ell^\infty$ input-output gain.

The signed identity \eqref{eq:signedmass} suggests measuring the memory-feedback amplification by the total variation of the Green sequence.  We therefore introduce the following quantity.

\begin{definition}\label{def:Lambda3}
For $\mu>0$, we set
\begin{equation}\label{eq:Lambda3mu}
 \LamthreeMu=\eta(\mu)\sum_{m=0}^{\infty}\abs{G_m(\mu)},
\end{equation}
and define the global BDF3 stability constant by
\begin{equation}\label{eq:Lambda3}
 \Lamthree=\sup_{\mu>0}\LamthreeMu.
\end{equation}
\end{definition}

\begin{lemma}\label{lem:Green-identities}

For every $\mu>0$, the Green sequence \eqref{eq:Greenrec} satisfies the generating-function identity \eqref{eq:GreenGF}; the series is absolutely convergent for $\abs{z}\le1$, and the signed-mass identity \eqref{eq:signedmass} holds.  Moreover,
\begin{equation}\label{eq:negmass}
 \LamthreeMu
 =1+2\eta(\mu)\sum_{\{m:\,G_m(\mu)<0\}}\abs{G_m(\mu)}.
\end{equation}
Finally, for every $\mu>3$, $\sum_{m=1}^{\infty}\abs{G_m(\mu)}\le \frac{29}{6\mu-18}$,
and hence
\[
 \sum_{m=1}^{\infty}\abs{G_m(\mu)}=O(\mu^{-1})
 \qquad (\mu\to\infty).
\]
\end{lemma}

\begin{proof}
Fix $\mu>0$ and suppress the dependence on $\mu$ in the notation.  Let
\[
 \mathcal{G}_{\mu}(z):=\sum_{m=0}^{\infty}G_m(\mu) z^m.
\]
Multiplying the recurrence in \eqref{eq:Greenrec} by $z^m$ and summing from $m=3$ to infinity gives
\begin{align*}
0={}&\sum_{m=3}^{\infty}G_mz^m
-a_1\sum_{m=3}^{\infty}G_{m-1}z^m
+a_2\sum_{m=3}^{\infty}G_{m-2}z^m
-a_3\sum_{m=3}^{\infty}G_{m-3}z^m\\
={}&\mathcal{G}_{\mu}(z)-G_0-G_1z-G_2z^2
-a_1z\bigl(\mathcal{G}_{\mu}(z)-G_0-G_1z\bigr)\\
&\quad +a_2z^2\bigl(\mathcal{G}_{\mu}(z)-G_0\bigr)-a_3z^3\mathcal{G}_{\mu}(z).
\end{align*}
Therefore, we get
\[
 \bigl(1-a_1z+a_2z^2-a_3z^3\bigr)\mathcal{G}_{\mu}(z)
 =G_0+(G_1-a_1G_0)z+(G_2-a_1G_1+a_2G_0)z^2.
\]
The initial values in \eqref{eq:Greenrec} make the coefficients of $z$ and $z^2$ vanish and give $G_0=1$.  This proves \eqref{eq:GreenGF} first as an identity of formal power series.  By Lemma~\ref{lem:schur}, all characteristic roots lie strictly in the unit disk, so the Green coefficients decay geometrically for each fixed $\mu>0$; consequently \eqref{eq:GreenGF} is absolutely convergent for $\abs{z}\le1$.

We may therefore set $z=1$ in \eqref{eq:GreenGF}.  Using \eqref{eq:coeffs}, we get
\[
 1-a_1(\mu)+a_2(\mu)-a_3(\mu)
 =\frac{d_\mu-11}{d_\mu}
 =\eta(\mu).
\]
Thus, we have
\[
 \sum_{m=0}^{\infty}G_m(\mu)=\frac1{\eta(\mu)},
\]
which is exactly \eqref{eq:signedmass}.

Next, we define the positive and negative masses
\[
 P_\mu:=\sum_{\{m:\,G_m(\mu)\ge0\}}G_m(\mu),
 \qquad
 N_\mu:=\sum_{\{m:\,G_m(\mu)<0\}}\abs{G_m(\mu)}.
\]
Then
\[
 \sum_{m\ge0}G_m=P_\mu-N_\mu,
 \qquad
 \sum_{m\ge0}\abs{G_m}=P_\mu+N_\mu.
\]
By \eqref{eq:signedmass}, $\eta(\mu)(P_\mu-N_\mu)=1$, and therefore
\[
 \LamthreeMu
 =\eta(\mu)(P_\mu+N_\mu)
 =\eta(\mu)(P_\mu-N_\mu)+2\eta(\mu)N_\mu
 =1+2\eta(\mu)N_\mu,
\]
which proves \eqref{eq:negmass}.

It remains to establish the large-$\mu$ estimate. We set
\[
 q_\mu(z):=a_1(\mu)z-a_2(\mu)z^2+a_3(\mu)z^3.
\]
Then \eqref{eq:GreenGF} reads $\mathcal{G}_{\mu}(z)=(1-q_\mu(z))^{-1}$.  In the Wiener algebra, with
\[
   \Big\|\sum_{j\ge0}c_jz^j \Big\|_{W}:=\sum_{j\ge0}\abs{c_j},
\]
one has $\|fg\|_W\le\|f\|_W\|g\|_W$. For $f(z)=\sum_{j\ge0}c_jz^j$ and $g(z)=\sum_{j\ge0}d_jz^j$,
the Cauchy product gives
\[
 \|fg\|_W\le\sum_{n\ge0}\sum_{j=0}^n|c_jd_{n-j}|
 =\|f\|_W\|g\|_W.
\]
  Since
\[
 \|q_\mu\|_W=a_1(\mu)+a_2(\mu)+a_3(\mu)
 =\frac{29}{11+6\mu}<1
 \qquad (\mu>3),
\]
the Neumann series converges in the Wiener algebra and yields
\begin{align*}
 \sum_{m=1}^{\infty}\abs{G_m(\mu)}
 &=\|\mathcal{G}_{\mu}-1\|_W
 = \Big\|\sum_{j=1}^{\infty}q_\mu^j \Big\|_W \leq \sum_{j=1}^{\infty}\big\|q_\mu^j \big\|_W   \\
 &\le \sum_{j=1}^{\infty}\|q_\mu\|_W^j
 =\frac{\|q_\mu\|_W}{1-\|q_\mu\|_W}
 =\frac{29}{6\mu-18},
\end{align*}
where we used $\mathcal{G}_{\mu}-1=\frac{1}{1-q_{\mu}} -1 = \sum_{j=1}^{\infty}q_{\mu}^j$ at the second equals sign. The stated $O(\mu^{-1})$ estimate follows immediately.
\end{proof}

By \eqref{eq:negmass}, $\LamthreeMu-1$ equals twice the weighted negative Green mass.  Hence $\LamthreeMu\ge1$, with equality precisely when all $G_m(\mu)$ are nonnegative.  The excess therefore measures oscillatory amplification.

\begin{proposition}\label{prop:LambdaFinite}
The function $\mu\mapsto\LamthreeMu$ is finite and continuous on $(0,\infty)$, extends continuously to the two endpoints after compactification, and satisfies
\begin{equation}\label{eq:LambdaLimits}
 \lim_{\mu\to0^+}\LamthreeMu=1,\qquad \lim_{\mu\to\infty}\LamthreeMu=1.
\end{equation}
In particular, the constant $\Lamthree$ in \eqref{eq:Lambda3} is finite and is attained at least at one interior point.
\end{proposition}

\begin{proof}
For every compact interval $[\mu_1,\mu_2]\subset(0,\infty)$, Lemma~\ref{lem:schur} and compactness give a uniform spectral radius strictly below one for the companion matrices associated with \eqref{eq:Greenrec}.  The roots depend continuously on $\mu$ and are simple, since the
discriminant $-108(36\mu^2-30\mu+13)$ is negative.  Their maximal
modulus is bounded by some $\rho_*<1$ on the compact interval; hence
$|G_m(\mu)|\le C\rho^m$ uniformly for any $\rho\in(\rho_*,1)$.
 Hence, the series in \eqref{eq:Lambda3mu} converges uniformly on that interval, which proves continuity there.

We next consider $\mu\to\infty$.  Lemma~\ref{lem:Green-identities} gives
\[
 \sum_{m=1}^{\infty}\abs{G_m(\mu)}=O(\mu^{-1}).
\]
Since $\eta(\mu)\to1$, the second limit in \eqref{eq:LambdaLimits} follows.

It remains to study the singular endpoint $\mu\to0^+$.  It is useful to keep the principal root and the two parasitic roots separate, since only the former approaches the unit circle.  Multiply the characteristic polynomial \eqref{eq:charpoly} by $d_\mu=11+6\mu$ and set
\[
 F(r,\mu):=(11+6\mu)r^3-18r^2+9r-2.
\]
At $\mu=0$, we have
\[
 F(r,0)=11r^3-18r^2+9r-2=(r-1)(11r^2-7r+2),
\]
so the principal root is $r_1(0)=1$, whereas the two parasitic roots are
\[
 r_{2,3}(0)=\frac{7\pm \mathrm{i}\sqrt{39}}{22},
 \qquad
 \abs{r_{2,3}(0)}=\sqrt{\frac{2}{11}}<1.
\]
The root $r=1$ is simple because
\[
 F_r(1,0)=33-36+9=6\neq0.
\]
Hence, by the analytic implicit-function theorem, there is a real analytic branch $r_1(\mu)$ for $\mu$ sufficiently small such that
\[
 F(r_1(\mu),\mu)=0,\qquad r_1(0)=1.
\]
Differentiating the identity $F(r_1(\mu),\mu)=0$ at $\mu=0$ gives
\[
 r_1'(0)
 =-\frac{F_\mu(1,0)}{F_r(1,0)}
 =-\frac{6}{6}=-1,
\]
which combines $r_1(0)=1$ to obtain
\begin{equation}\label{eq:principal-root-small-mu}
 r_1(\mu)=1-\mu+O(\mu^2)  \,\Rightarrow \,
 1-r_1(\mu)=\mu+O(\mu^2).
\end{equation}
In particular, $0<r_1(\mu)<1$ for all sufficiently small $\mu>0$.

The two parasitic roots are simple at $\mu=0$ and lie strictly inside the unit disk.  By continuity of the roots with respect to the polynomial coefficients, there exist $\mu_0>0$ and $\rho\in(0,1)$ such that
\begin{equation}\label{eq:parasitic-uniform-gap}
 \abs{r_2(\mu)},\ \abs{r_3(\mu)}\le \rho,
 \qquad 0\le\mu\le\mu_0.
\end{equation}
After reducing $\mu_0$ if necessary, the three roots remain distinct on $[0,\mu_0]$.  Vieta's formulas and partial fractions give
\[
 \mathcal G_\mu(z)=\prod_{j=1}^3(1-r_j(\mu)z)^{-1}
 =\sum_{j=1}^3\frac{A_j(\mu)}{1-r_j(\mu)z},\qquad
 A_j(\mu)=\frac{r_j(\mu)^2}{p_\mu'(r_j(\mu))}.
\]
Expanding each fraction as a geometric series yields
\begin{equation}\label{eq:Green-root-decomp}
 G_m(\mu)=\sum_{j=1}^3 A_j(\mu)r_j(\mu)^m.
\end{equation}
At $\mu=0$, the generating function has the factorization
\[
 \mathcal{G}_0(z)
 =\frac{1}{
 (1-z)\left(1-\frac{7}{11}z+\frac{2}{11}z^2\right)},
\]
so the coefficient of the principal pole is
\[
 A_1(0)
 =\lim_{z\to1}(1-z)\mathcal{G}_0(z)
 =\frac{1}{1-\frac{7}{11}+\frac{2}{11}}
 =\frac{11}{6}.
\]
Consequently, we have
\begin{equation}\label{eq:A1-small-mu}
 A_1(\mu)=\frac{11}{6}+O(\mu),
\end{equation}
while $A_2(\mu)$ and $A_3(\mu)$ remain uniformly bounded for $0\le\mu\le\mu_0$.

We can now estimate the absolute Green mass.  From \eqref{eq:Green-root-decomp}, \eqref{eq:parasitic-uniform-gap}, and $0<r_1(\mu)<1$, we have
\begin{align*}
 \sum_{m=0}^{\infty}\abs{G_m(\mu)}
 &\le
 \abs{A_1(\mu)}\sum_{m=0}^{\infty}r_1(\mu)^m
 +\sum_{j=2}^3\abs{A_j(\mu)}
   \sum_{m=0}^{\infty}\abs{r_j(\mu)}^m\le
 \frac{\abs{A_1(\mu)}}{1-r_1(\mu)}+C,
\end{align*}
where $C$ is independent of sufficiently small $\mu$. Therefore, we yield
\begin{equation}\label{eq:Lambda-small-mu-upper}
 \LamthreeMu
 \le
 \eta(\mu)\frac{\abs{A_1(\mu)}}{1-r_1(\mu)}
 +C\eta(\mu).
\end{equation}
Since
\[
 \eta(\mu)=\frac{6\mu}{11+6\mu}
 =\frac{6}{11}\mu+O(\mu^2),
\]
the relations \eqref{eq:principal-root-small-mu} and \eqref{eq:A1-small-mu} give
\[
 \eta(\mu)\frac{\abs{A_1(\mu)}}{1-r_1(\mu)}
 =
 \frac{
 \left(\frac{6}{11}\mu+O(\mu^2)\right)
 \left(\frac{11}{6}+O(\mu)\right)}
 {\mu+O(\mu^2)}
 \longrightarrow 1,
\]
whereas $C\eta(\mu)\to0$.  Hence, \eqref{eq:Lambda-small-mu-upper} implies
\[
 \limsup_{\mu\to0^+}\LamthreeMu\le1.
\]
On the other hand, the signed-mass identity \eqref{eq:signedmass} and the triangle inequality give
\[
 1
 =\left|\eta(\mu)\sum_{m=0}^{\infty}G_m(\mu)\right|
 \le \eta(\mu)\sum_{m=0}^{\infty}\abs{G_m(\mu)}
 =\LamthreeMu,
\]
so that
\[
 \liminf_{\mu\to0^+}\LamthreeMu\ge1.
\]
The first limit in \eqref{eq:LambdaLimits} follows.

It remains only to justify the final global statement.  The two endpoint limits show that $\LamthreeMu$ is bounded outside a compact subinterval of $(0,\infty)$, while continuity gives boundedness on every compact subinterval; hence $\Lamthree<\infty$.  If $\LamthreeMu\equiv1$, then the supremum is attained at every interior point.  Otherwise, since \eqref{eq:signedmass} implies $\LamthreeMu\ge1$, there exists $\mu_\ast>0$ with $\Lambda_3(\mu_\ast)>1$.  The endpoint limits allow us to choose $0<\delta<\mu_\ast<M<\infty$ such that
\[
 \LamthreeMu<\Lambda_3(\mu_\ast)
 \qquad\text{for }\mu\in(0,\delta]\cup[M,\infty).
\]
Therefore the global supremum is the maximum of the continuous function $\Lambda_3(\mu)$ on the compact interval $[\delta,M]$, and it is attained at an interior point of $(0,\infty)$.
\end{proof}

\begin{theorem}\label{thm:Lambda-enclosure}
The constant in Definition~\ref{def:Lambda3} satisfies
\begin{equation}\label{eq:Lambda-enclosure}
 \frac{23}{20}<\Lamthree<\frac{579}{500}.
\end{equation}
Consequently, $\beta<500/579$ is a fully explicit sufficient condition for the resolvent-based condition $\Lamthree\beta<1$.
\end{theorem}

\begin{proof}
For the lower bound, choose $\mu=15/2$, so that $\eta=45/56$.  Every finite partial sum of the absolute Green series gives
\[
 \Lamthree\ge\Lambda_3(15/2)
 \ge\frac{45}{56}\sum_{m=0}^{9}|G_m(15/2)|.
\]
At $\mu=15/2$, the recurrence gives
\[
    G_0=1,\qquad
    G_1=\frac{9}{28},\qquad
    G_2=
    \left(\frac{9}{28}\right)^2-\frac{9}{56}
    =-\frac{45}{784},
\]
and, for $m\ge3$,
\[
    G_m
    =
    \frac{9}{28}G_{m-1}
    -\frac{9}{56}G_{m-2}
    +\frac{1}{28}G_{m-3}.
\]
Successive application of this recurrence shows that
$G_0,\ldots,G_9$ have the sign pattern
$
    ++--++--++.
$
Consequently, the first ten terms can be evaluated exactly using
rational arithmetic, and one obtains
\[
    \frac{45}{56}
    \sum_{m=0}^{9}
    \left|G_m\!\left(\frac{15}{2}\right)\right|
    =
    \frac{681305421762045}{592393533390848} \approx 1.150089.
\]
Hence
\[
   \Lambda_3 \geq  \frac{681305421762045}{592393533390848}
    >
    \frac{23}{20}.
\]

The upper bound is established in Appendix~\ref{app:Lambda-bound},
where the parameter range is divided into suitable subintervals and
the corresponding Green-function contributions are bounded by exact
Bernstein enclosures together with analytic tail estimates.  These
estimates yield
\[
    \Lambda_3<\frac{579}{500}.
\]
Finally, $\beta<500/579$ and $\Lamthree<579/500$ imply $\Lamthree\beta<1$. We complete the proof.
\end{proof}

\begin{remark}\label{rem:LambdaNumerical}
Section~\ref{sec:numerics} evaluates \eqref{eq:Lambda3mu} with a controlled tail and finds
\[
 \Lamthree\approx1.1501238,
\]
with a maximizer near $\mu=7.4296$.  Hence the numerically indicated threshold is $\beta\approx0.86947$, whereas the rigorous analytic enclosure in Theorem~\ref{thm:Lambda-enclosure} yields the slightly more conservative but fully explicit condition $\beta<500/579\approx0.86356$.
\end{remark}

\subsection{Weighted Green bounds}

We next establish exponential resolvent bounds uniformly across the positive spectrum.

\begin{lemma}\label{lem:weighted-resolvent}
Let $\lambda_0>0$ and $\varepsilon>0$ be fixed. The constants below may depend on $\lambda_0$ and $\varepsilon$, but not on
$\lambda$, $k$, $n$, or an upper spectral cutoff.
  There exist $\gamma_0>0$, $k_0>0$, and $C>0$ such that, whenever $0<k\le k_0$, $0\le\gamma\le\gamma_0$, and $\mu=k\lambda$ with $\lambda\ge\lambda_0$, one has
\begin{equation}\label{eq:weighted-Lambda}
 \eta(\mu)\sum_{m=0}^{\infty}e^{\gamma km}\abs{G_m(\mu)}\le\Lamthree+\varepsilon
\end{equation}
and
\begin{equation}\label{eq:weighted-pointwise}
 \sup_{m\ge0}e^{\gamma km}\abs{G_m(\mu)}\le C.
\end{equation}
Moreover, every solution $h^n$ of the homogeneous recurrence associated with \eqref{eq:scalar-rec} satisfies
\begin{equation}\label{eq:homogeneous-bound}
 \abs{h^n}\le C e^{-\gamma t_n}\bigl(\abs{h^0}+\abs{h^1}+\abs{h^2}\bigr).
\end{equation}
\end{lemma}
\begin{proof}
  The proof, which requires a uniform discussion of the regimes $\mu\ll1$, $\mu=O(1)$, and $\mu\gg1$, is deferred to Appendix~\ref{app:weighted}.
\end{proof}

\section{Long-time stability}\label{sec:stability}

We combine the scalar resolvent bounds with exponential kernel decay.  For $\gamma\ge0$, we define
\[
 b_\gamma(t)=e^{\gamma t}\abs{b(t)},
\]
and let
\begin{equation}\label{eq:discrete-beta}
 \beta_{\gamma,k}^n:=\sigma^n\bigl(\{e^{\gamma(t_n-t_j)}\abs{b(t_n-t_j)}\}_{j=0}^{n-1}\bigr).
\end{equation}
Since the global weights in \eqref{eq:quadrature-global} are nonnegative, $\beta_{\gamma,k}^n$ is a genuine discrete weighted mass.

\begin{lemma}\label{lem:discrete-beta}
Let $b$ satisfy Assumption~\ref{ass:kernel} in its stability form and suppose \eqref{eq:kernel-sharp-condition}.  Then there exist numbers $\betabar$, $\gamma>0$, and $k_0>0$ such that
\begin{equation}\label{eq:beta-bar}
 \beta<\betabar<\Lamthree^{-1},\qquad
 \sup_{n\ge3}\beta_{\gamma,k}^n\le\betabar
\end{equation}
for every $0<k\le k_0$.
\end{lemma}

\begin{proof}
Condition \eqref{eq:kernel-sharp-condition} gives $\beta<\Lambda_3^{-1}$.  We may therefore choose a number
$\bar\beta$ such that
\[
    \beta<\bar\beta<\Lambda_3^{-1}.
\]
Set $\delta:=\frac{\bar\beta-\beta}{2}>0$. By Assumption~\ref{ass:kernel} in its stability form, there exist
$\alpha>0$ and constants $K_0,K_1>0$ such that
\[
    |b(t)|\le K_0e^{-\alpha t},
    \qquad
    |b'(t)|\le K_1e^{-\alpha t},
    \qquad t\ge0.
\]
For $0\le\gamma<\alpha$, we define $\beta_\gamma$ by \eqref{eq:beta-gamma}.
Since $e^{\gamma s}|b(s)|\to |b(s)|$ pointwise as $\gamma\to0^+$,
and, for example, for $0\le\gamma\le\alpha/2$,
\[
    e^{\gamma s}|b(s)|
    \le
    K_0e^{-\alpha s/2}\in L^1(0,\infty),
\]
the dominated convergence theorem yields
\[
    \beta_\gamma\longrightarrow\beta
    \qquad\text{as }\gamma\to0^+.
\]
Thus, we may choose $\gamma>0$, with $\gamma<\alpha/2$, so small that
\begin{equation}\label{eq:beta-gamma-margin}
    \beta_\gamma
    \le \beta+\delta
    =\frac{\beta+\bar\beta}{2}
    <\bar\beta.
\end{equation}

We compare the discrete weighted mass with $\beta_\gamma$.  For a
fixed $n\ge3$, define
\[
    g_n(s)
    :=
    e^{\gamma(t_n-s)}
    |b(t_n-s)|,
    \qquad 0\le s\le t_n.
\]
By \eqref{eq:discrete-beta}, we have $\beta_{\gamma,k}^n=\sigma^n(g_n)$. Since $b\in C^1([0,\infty))$, the function $|b|$ is absolutely continuous and
\[
    \frac{d}{dt}|b(t)|
    =
    \operatorname{sgn}(b(t))\,b'(t)
\]
for almost every $t\ge0$. Therefore, $g_n\in W^{1,1}(0,t_n)$ and, with
$\tau=t_n-s$,
\[
\begin{aligned}
    g_n'(s)
    &=
    -e^{\gamma\tau}
    \left[
        \gamma |b(\tau)|
        +
        \operatorname{sgn}(b(\tau))\,b'(\tau)
    \right]
\end{aligned}
\]
for almost every $s\in(0,t_n)$.  Consequently, we get
\[
    |g_n'(s)|
    \le
    e^{\gamma(t_n-s)}
    \left(
        \gamma |b(t_n-s)|
        +
        |b'(t_n-s)|
    \right).
\]
Applying Lemma~\ref{lem:quadrature} with $j=1$ gives
\[
    \left|
        \sigma^n(g_n)
        -
        \int_0^{t_n}g_n(s)\,ds
    \right|
    \le
    Ck\int_0^{t_n}|g_n'(s)|\,ds,
\]
where $C$ is independent of $k$ and $n$. Thus,
\begin{equation}\label{eq:discrete-beta-first}
\begin{aligned}
    \beta_{\gamma,k}^n
    &\le
    \int_0^{t_n}
        e^{\gamma(t_n-s)}
        |b(t_n-s)|\,ds
    \\
    &\quad
    +Ck
    \int_0^{t_n}
        e^{\gamma(t_n-s)}
        \left(
            \gamma|b(t_n-s)|
            +
            |b'(t_n-s)|
        \right)\,ds.
\end{aligned}
\end{equation}
Making the change of variables $\tau=t_n-s$ in both integrals, we obtain
\[
    \int_0^{t_n}
        e^{\gamma(t_n-s)}
        |b(t_n-s)|\,ds
    =
    \int_0^{t_n}e^{\gamma\tau}|b(\tau)|\,d\tau
    \le\beta_\gamma,
\]
and
\[
\begin{aligned}
    \int_0^{t_n}
        e^{\gamma(t_n-s)}
        \left(
            \gamma|b(t_n-s)|
            +
            |b'(t_n-s)|
        \right)\,ds
    &\le
    \int_0^\infty
        e^{\gamma\tau}
        \left(
            \gamma|b(\tau)|
            +
            |b'(\tau)|
        \right)\,d\tau.
\end{aligned}
\]
The last integral is finite.  Indeed, by the exponential fading
assumption and $\gamma<\alpha$,
\[
\begin{aligned}
    M_\gamma
    & :=
    \int_0^\infty
        e^{\gamma\tau}
        \left(
            \gamma|b(\tau)|
            +
            |b'(\tau)|
        \right)\,d\tau
    \\
    &\le \int_0^\infty
        \left(\gamma K_0+K_1\right)
        e^{-(\alpha-\gamma)\tau}\,d\tau
    =
    \frac{\gamma K_0+K_1}{\alpha-\gamma}
    <\infty.
\end{aligned}
\]
Substituting these bounds into \eqref{eq:discrete-beta-first} gives
\begin{equation}\label{eq:discrete-beta-uniform}
    \beta_{\gamma,k}^n
    \le
    \beta_\gamma
    +
    CkM_\gamma,
    \qquad n\ge3.
\end{equation}
Finally, choose $k_0>0$ sufficiently small that $Ck_0M_\gamma\le\delta$. Then, for every $0<k\le k_0$ and every $n\ge3$, combining
\eqref{eq:beta-gamma-margin} and \eqref{eq:discrete-beta-uniform} yields
\[
\begin{aligned}
    \beta_{\gamma,k}^n
    &\le
    \beta_\gamma+CkM_\gamma
    \le
    (\beta+\delta)+\delta
    =
    \bar\beta,
\end{aligned}
\]
which, together with the initial choice $\beta<\bar\beta<\Lambda_3^{-1}$, proves \eqref{eq:beta-bar}.
\end{proof}

\begin{theorem}\label{thm:stability}
Let $b$ satisfy Assumption~\ref{ass:kernel} in its stability form and
assume \eqref{eq:kernel-sharp-condition}.  Then there exist constants
$C>0$, $\gamma>0$, and $k_0>0$, independent of the final time, such
that the solution of \eqref{eq:scheme} satisfies
\begin{equation}\label{eq:stability}
 \norm{U^n}
 \le
 C e^{-\gamma t_n}
 \bigl(\norm{U^0}+\norm{U^1}+\norm{U^2}\bigr)
 +
 Ck\sum_{j=3}^n
 e^{-\gamma(t_n-t_j)}\norm{f^j},
\end{equation}
for all $n\ge3$ and $0<k\le k_0$.
\end{theorem}

\begin{proof}
We choose $\bar\beta$ as in Lemma~\ref{lem:discrete-beta}, so that
\[
 \beta<\bar\beta<\Lambda_3^{-1},
 \qquad
 \sup_{n\ge3}\beta_{\gamma,k}^n\le\bar\beta
\]
for all sufficiently small $k$.  Since
$\Lambda_3\bar\beta<1$, we may choose $\varepsilon>0$ such that
\begin{equation}\label{eq:absorption}
 (\Lambda_3+\varepsilon)\bar\beta<1.
\end{equation}
By Lemma~\ref{lem:weighted-resolvent}, after possibly decreasing
$\gamma$ and $k_0$, the estimates
\eqref{eq:weighted-Lambda}--\eqref{eq:homogeneous-bound} hold with
this value of $\varepsilon$.

We first work with a fixed eigenmode.  Let
$c_m^n=\inner{U^n}{\varphi_m}$,
$f_m^n=\inner{f^n}{\varphi_m}$, and
$\mu_m=k\lambda_m$.  The scalar recurrence associated with
\eqref{eq:scheme} takes the form
\[
 c_m^n-a_1(\mu_m)c_m^{n-1}
 +a_2(\mu_m)c_m^{n-2}
 -a_3(\mu_m)c_m^{n-3}
 =
 \eta(\mu_m)\sigma^n(b_nc_m)
 +\frac{6k}{d_{\mu_m}}f_m^n .
\]
We set
\[
 \psi_m^j
 :=
 \eta(\mu_m)\sigma^j(b_jc_m)
 +\frac{6k}{d_{\mu_m}}f_m^j,
\]
and let $h_m^n$ denote the solution of the corresponding homogeneous
recurrence with $h_m^\ell=c_m^\ell$, $\ell=0,1,2$. To obtain the discrete variation-of-constants formula, put
$w_m^n:=c_m^n-h_m^n$.  Then $w_m^0=w_m^1=w_m^2=0$ and
\[
 w_m^n-a_1w_m^{n-1}
 +a_2w_m^{n-2}
 -a_3w_m^{n-3}
 =
 \psi_m^n .
\]
The Green sequence satisfies
\[
 G_0=1,\qquad
 G_1-a_1G_0=0,\qquad
 G_2-a_1G_1+a_2G_0=0,
\]
and, for $\ell\ge3$,
\[
 G_\ell-a_1G_{\ell-1}
 +a_2G_{\ell-2}
 -a_3G_{\ell-3}=0.
\]
Hence, direct substitutions show that
\[
 w_m^n
 =
 \sum_{j=3}^n
 G_{n-j}(\mu_m)\psi_m^j.
\]
Consequently, we get
\begin{equation}\label{eq:variation}
\begin{aligned}
 c_m^n
 ={}&
 h_m^n
 +\eta(\mu_m)
 \sum_{j=3}^n
 G_{n-j}(\mu_m)\sigma^j(b_jc_m)
 +\frac{6k}{d_{\mu_m}}
 \sum_{j=3}^n
 G_{n-j}(\mu_m)f_m^j .
\end{aligned}
\end{equation}
Define the exponentially weighted modal maximum
\[
 \mathcal M_m^n
 :=
 \max_{0\le\ell\le n}
 e^{\gamma t_\ell}|c_m^\ell|.
\]
We first estimate the memory term.  Since the assembled quadrature
weights are nonnegative and the quadrature at level $j$ uses only
values $c_m^\ell$ with $\ell\le j-1$, we have
\begin{equation}\label{eq:memory-mode}
\begin{aligned}
 e^{\gamma t_j}
 \abs{\sigma^j(b_jc_m)}
 &\le
 \sigma^j
 \left(
 \left\{
 e^{\gamma(t_j-t_\ell)}
 |b(t_j-t_\ell)|
 \,e^{\gamma t_\ell}|c_m^\ell|
 \right\}_{\ell=0}^{j-1}
 \right)
 \\
 &\le
 \beta_{\gamma,k}^j \mathcal M_m^{j-1}
 \le
 \bar\beta \mathcal M_m^n.
\end{aligned}
\end{equation}

Now fix an arbitrary $r$ with $3\le r\le n$.  We then multiply
\eqref{eq:variation} by $e^{\gamma t_r}$ and use
$t_r-t_j=(r-j)k$ to arrive at
\begin{align*}
 e^{\gamma t_r}|c_m^r|
 \le{}&
 e^{\gamma t_r}|h_m^r|
 +
 \eta(\mu_m)
 \sum_{j=3}^r
 e^{\gamma k(r-j)}
 |G_{r-j}(\mu_m)|
 e^{\gamma t_j}
 |\sigma^j(b_jc_m)|
 \\
 &+
 \frac{6k}{d_{\mu_m}}
 \sum_{j=3}^r
 e^{\gamma k(r-j)}
 |G_{r-j}(\mu_m)|
 e^{\gamma t_j}|f_m^j|.
\end{align*}
For the homogeneous part,
\eqref{eq:homogeneous-bound} yields
\[
 e^{\gamma t_r}|h_m^r|
 \le
 C\sum_{\ell=0}^2|c_m^\ell|.
\]
For the memory part, \eqref{eq:memory-mode} and
\eqref{eq:weighted-Lambda} give
\begin{align*}
 &\eta(\mu_m)
 \sum_{j=3}^r
 e^{\gamma k(r-j)}
 |G_{r-j}(\mu_m)|
 \Big(e^{\gamma t_j}
 |\sigma^j(b_jc_m)| \Big)
 \\
 &\qquad\le
 \bar\beta \mathcal M_m^n\,
 \eta(\mu_m)
 \sum_{j=3}^r
 e^{\gamma k(r-j)}
 |G_{r-j}(\mu_m)|
 \le
 (\Lambda_3+\varepsilon)\bar\beta \mathcal M_m^n.
\end{align*}
Finally, by \eqref{eq:weighted-pointwise}, we have $\sup_{\ell\ge0}
 e^{\gamma k\ell}|G_\ell(\mu_m)| \le C$, and, since $d_{\mu_m}=11+6\mu_m\ge11$, we obtain
\begin{align*}
 &\frac{6k}{d_{\mu_m}}
 \sum_{j=3}^r
 e^{\gamma k(r-j)}
 |G_{r-j}(\mu_m)|
 e^{\gamma t_j}|f_m^j|
 \\
 &\qquad\qquad\le
 Ck\sum_{j=3}^r e^{\gamma t_j}|f_m^j|
 \le
 Ck\sum_{j=3}^n e^{\gamma t_j}|f_m^j|.
\end{align*}
Combining these estimates gives, for every $3\le r\le n$, we get
\begin{equation}\label{eq:modal-premax}
 e^{\gamma t_r}|c_m^r|
 \le
 C\sum_{\ell=0}^2|c_m^\ell|
 +
 (\Lambda_3+\varepsilon)\bar\beta \mathcal M_m^n
 +
 Ck\sum_{j=3}^n e^{\gamma t_j}|f_m^j|.
\end{equation}
For $r=0,1,2$, the definition of $\mathcal M_m^n$ gives
\[
 e^{\gamma t_r}|c_m^r|
 \le
 e^{2\gamma k_0}
 \sum_{\ell=0}^2|c_m^\ell| \le C\sum_{\ell=0}^2|c_m^\ell|,
\]
We take the maximum over $0\le r\le n$ in
\eqref{eq:modal-premax} to yield
\[
 \mathcal M_m^n
 \le
 C\sum_{\ell=0}^2|c_m^\ell|
 +
 (\Lambda_3+\varepsilon)\bar\beta \mathcal M_m^n
 +
 Ck\sum_{j=3}^n e^{\gamma t_j}|f_m^j|.
\]
By \eqref{eq:absorption}, we give
\begin{equation}\label{eq:modal-stability}
 e^{\gamma t_n}|c_m^n|
 \le \mathcal M_m^n
 \le
 C\sum_{\ell=0}^2|c_m^\ell|
 +
 Ck\sum_{j=3}^n e^{\gamma t_j}|f_m^j|.
\end{equation}
We now return to the Hilbert-space norm. Taking the $\ell^2$ norm in
$m$ on both sides of \eqref{eq:modal-stability} and using the
triangle inequality in $\ell^2$, we have
\[
\begin{aligned}
 e^{\gamma t_n}
 \left(\sum_m|c_m^n|^2\right)^{1/2}
 \le{}&
 C\sum_{\ell=0}^2
 \left(\sum_m|c_m^\ell|^2\right)^{1/2}
 +
 Ck\sum_{j=3}^n
 e^{\gamma t_j}
 \left(\sum_m|f_m^j|^2\right)^{1/2}.
\end{aligned}
\]
By Parseval's identity, we have
\[
 \left(\sum_m|c_m^\ell|^2\right)^{1/2}=\|U^\ell\|,
 \qquad
 \left(\sum_m|f_m^j|^2\right)^{1/2}=\|f^j\|.
\]
Thus, we obtain
\[
 e^{\gamma t_n}\|U^n\|
 \le
 C\bigl(\|U^0\|+\|U^1\|+\|U^2\|\bigr)
 +
 Ck\sum_{j=3}^n e^{\gamma t_j}\|f^j\|.
\]
Multiplying by $e^{-\gamma t_n}$ yields
\[
 \|U^n\|
 \le
 Ce^{-\gamma t_n}
 \bigl(\|U^0\|+\|U^1\|+\|U^2\|\bigr)
 +
 Ck\sum_{j=3}^n
 e^{-\gamma(t_n-t_j)}\|f^j\|,
\]
which yields \eqref{eq:stability}. We thus complete the proof.
\end{proof}

\begin{corollary}\label{cor:explicit-stability}
Under the regularity assumptions of Theorem~\ref{thm:stability}, the estimate \eqref{eq:stability} holds whenever
\[
 \int_0^\infty\abs{b(s)}\,ds<\frac{500}{579}.
\]
\end{corollary}

\begin{proof}
Theorem~\ref{thm:Lambda-enclosure} gives $\Lamthree<579/500$, so the displayed condition implies \eqref{eq:kernel-sharp-condition}.
\end{proof}

\begin{remark}
No upper bound on the spectrum of $A$ appears in Theorem~\ref{thm:stability}.  Since the spectrum of $A$ is unbounded, the estimate involves no restriction of the form $k\lambda_{\max}\le C$.
 Thus the BDF3 time step does not have to satisfy a parabolic restriction of the form $k\lambda_{\max}\le C$.  This is consistent with the fact that the negative real axis belongs to the BDF3 stability region, but the proof is stronger than a standard absolute-stability argument because it also controls the memory feedback uniformly through the $\ell^1$ norm of the discrete Green function. The theorem still requires $0<k\le k_0$, as imposed by the memory-quadrature bounds.
\end{remark}

\section{Third-order convergence}\label{sec:convergence}

We next derive the temporal error estimate.  We begin by showing that the starting procedure \eqref{eq:start1}--\eqref{eq:start2} generates the accuracy required by BDF3.

\subsection{Accuracy of the starting values}

Let
\[
 M(t)=\int_0^t b(t-s)Au(s)\,ds.
\]
We assume in this subsection that the exact solution is smooth on $[0,2k]$ for the expressions below to be meaningful.

\begin{lemma}\label{lem:startup}
Let $U^0=u_0$, and let $U^1,U^2$ be defined by
\eqref{eq:start1} and \eqref{eq:start2}.  Suppose that, on some fixed interval $[0,T_0]$,
$u\in C^3([0,T_0];H)$, $Au\in C^2([0,T_0];H)$, and
$b\in C^2([0,T_0])$.  Then, for $0<2k\le T_0$ sufficiently small,
\begin{equation}\label{eq:start-error}
 \norm{U^1-u(t_1)}+\norm{U^2-u(t_2)}
 \le Ck^3\calS_0,
\end{equation}
where $\calS_0$ depends on the corresponding sup norms on the fixed
interval $[0,T_0]$ and is independent of $k$.
\end{lemma}

\begin{proof}
Set $e^j:=U^j-u(t_j)$ for $j=0,1,2$.
Since $U^0=u_0=u(0)$, we have $e^0=0$. We consider the step from $t_0$ to $t_1$.  Averaging the
continuous equations
\[
 u_t(t_j)+Au(t_j)=M(t_j)+f(t_j),
 \qquad j=0,1,
\]
and using $M(0)=0$ gives
\begin{equation}\label{eq:start-exact-average}
 \frac{u_t(t_1)+u_t(t_0)}{2}
 +\frac12A\bigl(u(t_1)+u(t_0)\bigr)
 =
 \frac12M(t_1)+\frac12(f^1+f^0).
\end{equation}
The centered trapezoidal approximation of the derivative satisfies
\[
 \left\|
 D_1u(t_1)-\frac{u_t(t_1)+u_t(t_0)}{2}
 \right\|
 \le
 Ck^2\max_{0\le t\le k}\norm{u^{(3)}(t)}.
\]
Indeed, since $D_1u(t_1) =\frac1k\int_0^k u_t(s)\,ds$, this is precisely the trapezoidal-rule error applied to $u_t$. For the memory term, we define $\phi_1(s):=b(k-s)Au(s)$ on $[0,k]$.  Then
$M(t_1)=\int_0^k\phi_1(s)\,ds$, while
$\phi_1(0)=b(k)Au(0)$.  Hence,
\[
\begin{aligned}
 \norm{M(t_1)-k\,b(k)Au(0)}
 &=
 \left\|
 \int_0^k\bigl(\phi_1(s)-\phi_1(0)\bigr)\,ds
 \right\|                                                  \\
 &\le
 \int_0^k\int_0^s\norm{\phi_1'(\xi)}\,d\xi\,ds
 \le Ck^2\calS_0,
\end{aligned}
\]
where
\[
 \phi_1'(s)
 =
 -b'(k-s)Au(s)+b(k-s)Au_t(s)
\]
is uniformly bounded on $[0,k]$ under the stated assumptions. Consequently, \eqref{eq:start-exact-average} can be rewritten in the
form
\[
 D_1u(t_1)
 +\frac12A\bigl(u(t_1)+u(t_0)\bigr)
 =
 \frac{k}{2}b(k)Au(0)
 +\frac12(f^1+f^0)
 +\tau_1,
\]
where
\begin{equation}\label{eq:start-tau1}
 \norm{\tau_1}\le Ck^2\calS_0.
\end{equation}
Subtracting the starting scheme \eqref{eq:start1} and using
$e^0=0$ gives $\left(k^{-1}I+\frac12A\right)e^1=-\tau_1$. Since $A$ is selfadjoint and positive, we have
\[
 \left\|
 \left(k^{-1}I+\frac12A\right)^{-1}
 \right\|
 \le k.
\]
Therefore, by \eqref{eq:start-tau1}, we get
\begin{equation}\label{eq:e1-start}
 \norm{e^1}
 \le
 k\norm{\tau_1}
 \le Ck^3\calS_0.
\end{equation}
We next consider the second starting step.  At $t_2=2k$, the continuous equation is
\[
 u_t(t_2)+Au(t_2)=M(t_2)+f^2.
\]
The BDF2 consistency error $\tau_2^{\rm BDF}:=D_2u(t_2)-u_t(t_2)$
satisfies that
\begin{equation}\label{eq:bdf2-start-error}
 \norm{\tau_2^{\rm BDF}}
 \le
 Ck^2\max_{0\le t\le2k}\norm{u^{(3)}(t)}.
\end{equation}
For the memory term, we set
\[
 \phi_2(s):=b(2k-s)Au(s),\qquad 0\le s\le2k.
\]
Then we have
\[
 M(t_2)=\int_0^{2k}\phi_2(s)\,ds,
 \qquad
 \phi_2(k)=b(k)Au(t_1).
\]
The midpoint rule on $[0,2k]$ thus gives
\begin{equation}\label{eq:memory-midpoint-start}
 M(t_2)
 =
 2k\,b(k)Au(t_1)+\tau_2^{\rm mem},
 \qquad
 \norm{\tau_2^{\rm mem}}\le Ck^3\calS_0.
\end{equation}
Indeed,
\[
 \phi_2''(s)
 =
 b''(2k-s)Au(s)
 -2b'(2k-s)Au_t(s)
 +b(2k-s)Au_{tt}(s),
\]
which is uniformly bounded on $[0,2k]$, and the midpoint-rule error is
of order $k^3$. Combining \eqref{eq:bdf2-start-error} and
\eqref{eq:memory-midpoint-start}, the exact solution satisfies that
\[
 D_2u(t_2)+Au(t_2)
 =
 2k\,b(k)Au(t_1)+f^2+\tau_2,
\]
where
\begin{equation}\label{eq:start-tau2}
 \norm{\tau_2}
 \le Ck^2\calS_0.
\end{equation}
Here the $O(k^3)$ memory error is absorbed into the $O(k^2)$ BDF2 error. Subtracting this relation from \eqref{eq:start2} yields
\[
 D_2e^2+Ae^2
 =
 2k\,b(k)Ae^1-\tau_2.
\]
Since
\[
 D_2e^2
 =
 \frac{3e^2-4e^1+e^0}{2k}
 =
 \frac{3}{2k}e^2-\frac{2}{k}e^1,
\]
we obtain the explicit error equation
\begin{equation}\label{eq:e2-start-equation}
 \left(\frac{3}{2k}I+A\right)e^2
 =
 \frac{2}{k}e^1
 +2k\,b(k)Ae^1
 -\tau_2.
\end{equation}
Let $R_k:=\left(\frac{3}{2k}I+A\right)^{-1}$, then the positivity of $A$ gives
\[
 \norm{R_k}\le\frac{2k}{3},
 \qquad
 \norm{R_kA}\le1.
\]
Applying $R_k$ to \eqref{eq:e2-start-equation} therefore yields
\[
\begin{aligned}
 \norm{e^2}
 &\le
 \frac{2k}{3}\frac{2}{k}\norm{e^1}
 +2k|b(k)|\,\norm{R_kAe^1}
 +\frac{2k}{3}\norm{\tau_2}                                      \\
 &\le
 \frac43\norm{e^1}
 +2k|b(k)|\,\norm{e^1}
 +Ck^3\calS_0.
\end{aligned}
\]
Since $b$ is bounded on the initial interval and
\eqref{eq:e1-start} gives $\norm{e^1}\le Ck^3\calS_0$, we conclude that
\[
 \norm{e^2}
 \le
 Ck^3\calS_0+Ck^4\calS_0
 \le
 Ck^3\calS_0
\]
for $k$ sufficiently small.  Combining this estimate with
\eqref{eq:e1-start} proves \eqref{eq:start-error}.
\end{proof}

\subsection{The consistency residual}

Let $e^n=U^n-u(t_n)$ with $n\ge3$. Subtracting the continuous equation from \eqref{eq:scheme}, we obtain
\begin{equation}\label{eq:error-equation}
 \bdfthree e^n+Ae^n=\sigma^n(b_nAe)-\tau^n+q^n(g_n),\qquad n\ge3,
\end{equation}
where
\begin{equation}\label{eq:tau}
 \tau^n=\bdfthree u(t_n)-u_t(t_n),\qquad g_n(s)=b(t_n-s)Au(s).
\end{equation}
The definitions in \eqref{eq:tau} separate the BDF3 consistency error from the memory quadrature error.  A Taylor expansion of the BDF3 formula \eqref{eq:BDF3} shows that the leading local error is $-k^3u^{(4)}(t_n)/4$.  For the long-time summation it is more convenient to use the integral estimate
\begin{equation}\label{eq:tau-bound}
 \norm{\tau^n}\le Ck^2\int_{t_{n-3}}^{t_n}\norm{u^{(4)}(s)}\,ds.
\end{equation}
The third derivative of $g_n$ is
\begin{equation}\label{eq:gthird}
\begin{aligned}
 \partial_s^3g_n(s)={}&-b^{(3)}(t_n-s)Au(s)+3b''(t_n-s)Au_t(s)\\
 &-3b'(t_n-s)Au_{tt}(s)+b(t_n-s)Au_{ttt}(s).
\end{aligned}
\end{equation}
Combining \eqref{eq:gthird} with Assumption~\ref{ass:kernel} and Lemma~\ref{lem:quadrature} gives
\begin{equation}\label{eq:qmemory-bound}
 \norm{q^n(g_n)}\le Ck^3\int_0^{t_n}e^{-\alpha(t_n-s)}\sum_{\ell=0}^3\norm{Au^{(\ell)}(s)}\,ds.
\end{equation}
For compact notation, we define
\begin{equation}\label{eq:R3}
 \calR_3(s)=\norm{u^{(4)}(s)}+\sum_{\ell=0}^3\norm{Au^{(\ell)}(s)}.
\end{equation}

\begin{theorem}\label{thm:convergence}
Assume the hypotheses of Theorem~\ref{thm:stability}, the kernel derivative bounds in Assumption~\ref{ass:kernel} for $0\le j\le3$, and enough regularity for \eqref{eq:R3} to be locally integrable.  Let $U^0=u_0$, let the starting errors satisfy \eqref{eq:start-error}, and let $U^n$ be defined by \eqref{eq:scheme}.  Then there are constants $C>0$, $\delta>0$, and $k_0>0$, independent of $n$ and the final time, such that
\begin{equation}\label{eq:main-error}
 \norm{U^n-u(t_n)}\le C e^{-\delta t_n}\bigl(\norm{e^1}+\norm{e^2}\bigr)
 +Ck^3\int_0^{t_n}e^{-\delta(t_n-s)}\calR_3(s)\,ds
\end{equation}
for every $n\ge3$ and $0<k\le k_0$.
\end{theorem}

\begin{proof}
Apply Theorem~\ref{thm:stability} to the error equation \eqref{eq:error-equation}.  We obtain
\begin{equation}\label{eq:error-stabstep}
 \norm{e^n}\le C e^{-\gamma t_n}\bigl(\norm{e^1}+\norm{e^2}\bigr)
 +Ck\sum_{j=3}^n e^{-\gamma(t_n-t_j)}\left(\norm{\tau^j}+\norm{q^j(g_j)}\right).
\end{equation}
Using \eqref{eq:tau-bound}, finite overlap of the intervals $[t_{j-3},t_j]$, and a harmless reduction of the exponential rate, we obtain
\[
 k\sum_{j=3}^n e^{-\gamma(t_n-t_j)}\norm{\tau^j}
 \le Ck^3\int_0^{t_n}e^{-\delta(t_n-s)}\norm{u^{(4)}(s)}\,ds.
\]
For the memory quadrature term, insert \eqref{eq:qmemory-bound} into \eqref{eq:error-stabstep}.  Changing the order of the discrete summation and the continuous integral produces the convolution of two exponentially decaying kernels.  If $0<\delta<\min\{\alpha,\gamma\}$, then
\[
 k\sum_{j=3}^n e^{-\gamma(t_n-t_j)}\int_0^{t_j}e^{-\alpha(t_j-s)}\Phi(s)\,ds
 \le C\int_0^{t_n}e^{-\delta(t_n-s)}\Phi(s)\,ds
\]
for every nonnegative $\Phi$.  Taking $\Phi(s)=\sum_{\ell=0}^3\norm{Au^{(\ell)}(s)}$ and combining the two estimates proves \eqref{eq:main-error}.
\end{proof}

\begin{corollary}\label{cor:uniform-third}
Under the assumptions of Theorem~\ref{thm:convergence}, suppose in addition that $\sup_{t\ge0}\calR_3(t)<\infty$.
Then
\begin{equation}\label{eq:uniform-third}
 \sup_{n\ge0}\norm{U^n-u(t_n)}\le Ck^3,
\end{equation}
where $C$ is independent of the final time.
\end{corollary}

\begin{proof}
Lemma~\ref{lem:startup} gives the required $O(k^3)$ starting errors.  The exponential convolution in \eqref{eq:main-error} is uniformly bounded, so \eqref{eq:uniform-third} follows immediately.
\end{proof}

\section{Finite element approximation}\label{sec:fe}

We apply the resolvent estimates to a conforming Galerkin method.  Let $\Omega\subset\mathbb{R}^d$ be a bounded $C^{1,1}$ domain (or a convex polygonal/polyhedral domain with the required elliptic regularity), let $A$ be generated by a symmetric coercive bilinear form $a(\cdot,\cdot)$ on $H_0^1(\Omega)$, and let $S_h\subset H_0^1(\Omega)$ be a conforming finite element space.  We define $A_h:S_h\to S_h$ and the $L^2$ projection $P_h$ by
\[
 \inner{A_hv_h}{\chi_h}=a(v_h,\chi_h),\qquad
 \inner{P_hv}{\chi_h}=\inner{v}{\chi_h},\qquad \chi_h\in S_h.
\]
The Ritz projection $R_h$ is characterized by
\[
 a(R_hv-v,\chi_h)=0\qquad\text{for all }\chi_h\in S_h.
\]
We assume the standard Ritz estimate \cite{Thomee2006},
\begin{equation}\label{eq:Ritz}
 \norm{R_hv-v}\le Ch^r\norm{v}_{H^r(\Omega)}
\end{equation}
for the regularity index used below.  Since $P_h$ is the identity on $S_h$, we shall repeatedly use $P_hR_h=R_h$ in addition to the familiar identity $A_hR_hv=P_hAv$.

The fully discrete BDF3 method is
\begin{equation}\label{eq:FE-scheme}
 \bdfthree U_h^n+A_hU_h^n=\sigma^n(b_nA_hU_h)+P_hf^n,
\end{equation}
with finite element analogues of the starting steps \eqref{eq:start1}--\eqref{eq:start2}.  Coercivity and the Poincar\'e inequality imply that the smallest eigenvalue of $A_h$ is bounded below by a positive constant independent of $h$.  Consequently, the weighted scalar resolvent estimates established above apply to $A_h$ with constants that are uniform in the spatial mesh.

For the error splitting, write
\[
 U_h^n-u(t_n)=\theta^n+\rho^n,
 \qquad \theta^n=U_h^n-R_hu(t_n),\qquad \rho^n=R_hu(t_n)-u(t_n).
\]
We derive the consistency equation for $\theta^n$.  Let
\begin{align*}
 \tau^n&=\bdfthree u(t_n)-u_t(t_n),\qquad
 g_n(s)=b(t_n-s)Au(s),\\
 q^n(g_n)&=\sigma^n(g_n)-\int_0^{t_n}g_n(s)\,ds.
\end{align*}
By the Ritz relation, $\sigma^n(b_nA_hR_hu)=P_h\sigma^n(b_nAu)$.  Moreover, we have
\[
 P_hu_t(t_n)-\bdfthree R_hu(t_n)
 =-P_h\tau^n-P_h\bdfthree\rho^n,
\]
because $P_hR_h=R_h$.  Subtracting the projected continuous equation from \eqref{eq:FE-scheme} therefore gives the exact identity
\begin{equation}\label{eq:FE-theta-equation}
 \bdfthree\theta^n+A_h\theta^n
 =\sigma^n(b_nA_h\theta)-P_h\tau^n-P_h\bdfthree\rho^n+P_hq^n(g_n),
 \qquad n\ge3.
\end{equation}
This representation is useful because the discrete derivative of the Ritz error can be estimated directly by time integrals of $u_t$, without imposing an additional pointwise-in-time approximation assumption on $u_{tt}$.

Indeed, the BDF3 operator can be expressed through three backward first differences as
\[
 \bdfthree v^n=\frac{11}{6}\bar\partial v^n-\frac{7}{6}\bar\partial v^{n-1}+\frac13\bar\partial v^{n-2},
 \qquad \bar\partial v^j=\frac{v^j-v^{j-1}}{k}.
\]
Since $\rho(t)=(R_h-I)u(t)$, the Ritz estimate \eqref{eq:Ritz} yields
\[
 \norm{\bar\partial\rho^j}
 \le \frac{Ch^r}{k}\int_{t_{j-1}}^{t_j}\norm{u_t(s)}_{H^r(\Omega)}\,ds,
\]
and hence
\begin{equation}\label{eq:D3-rho}
 \norm{\bdfthree\rho^j}
 \le \frac{Ch^r}{k}\int_{t_{j-3}}^{t_j}\norm{u_t(s)}_{H^r(\Omega)}\,ds,
 \qquad j\ge3.
\end{equation}
The factor $k^{-1}$ in \eqref{eq:D3-rho} is exactly compensated by the factor $k$ in the discrete variation-of-constants estimate; the overlap of the intervals $[t_{j-3},t_j]$ is uniformly bounded.

\begin{theorem}\label{thm:FE}
Assume the hypotheses of Theorem~\ref{thm:convergence}, the Ritz estimate \eqref{eq:Ritz}, and starting values satisfying that
\[
 \mathcal{E}_{\mathrm{start}}
 :=\sum_{j=0}^2\norm{U_h^j-R_hu(t_j)}<\infty.
\]
Then there exist $C>0$, $\delta>0$, $h_0>0$, and $k_0>0$, independent of $n$, $h$, and the final time, such that
\begin{align}\label{eq:FE-error}
 \norm{U_h^n-u(t_n)}
 &\le Ch^r\left(\norm{u(t_n)}_{H^r(\Omega)}+
 \int_0^{t_n}e^{-\delta(t_n-s)}\norm{u_t(s)}_{H^r(\Omega)}\,ds\right)\notag\\
 &\quad+Ck^3\int_0^{t_n}e^{-\delta(t_n-s)}\calR_3(s)\,ds
 +Ce^{-\delta t_n}\mathcal{E}_{\mathrm{start}},
\end{align}
for all $n\ge3$, $0<h\le h_0$, and $0<k\le k_0$.
\end{theorem}

\begin{proof}
The lower spectral bound for $A_h$ is uniform in $h$, so Theorem~\ref{thm:stability} applied to \eqref{eq:FE-theta-equation} gives
\begin{align*}
 \norm{\theta^n}
 &\le Ce^{-\delta t_n}\mathcal{E}_{\mathrm{start}}
 +Ck\sum_{j=3}^ne^{-\delta(t_n-t_j)}
 \bigl(\norm{\tau^j}+\norm{q^j(g_j)}+\norm{\bdfthree\rho^j}\bigr).
\end{align*}
The truncation and quadrature terms are bounded exactly as in the proof of Theorem~\ref{thm:convergence}.  For the last term, \eqref{eq:D3-rho} gives
\begin{align*}
 k\sum_{j=3}^ne^{-\delta(t_n-t_j)}\norm{\bdfthree\rho^j}
 &\le Ch^r\sum_{j=3}^ne^{-\delta(t_n-t_j)}
 \int_{t_{j-3}}^{t_j}\norm{u_t(s)}_{H^r}\,ds\\
 &\le Ch^r\int_0^{t_n}e^{-\delta_1(t_n-s)}\norm{u_t(s)}_{H^r}\,ds,
\end{align*}
where $0<\delta_1<\delta$ is fixed.  In the last step we used the uniformly bounded overlap of the three-step intervals and the elementary comparison $e^{-\delta(t_n-t_j)}\le C e^{-\delta_1(t_n-s)}$ for $s\in[t_{j-3},t_j]$ and $k\le k_0$.  Relabeling $\delta_1$ as $\delta$, and using \eqref{eq:Ritz} to bound
\[
 \norm{\rho^n}\le Ch^r\norm{u(t_n)}_{H^r(\Omega)},
\]
proves \eqref{eq:FE-error}.
\end{proof}

\begin{lemma}\label{lem:FE-start}
Assume the hypotheses of Lemma~\ref{lem:startup}, the Ritz estimate
\eqref{eq:Ritz}, and $u_t\in C([0,T_0];H^r(\Omega))$.
Take $U_h^0=R_hu_0$ and use the $A_h,P_h$ versions of
\eqref{eq:start1}--\eqref{eq:start2}.  Then for $2k\le T_0$,
\[
 \mathcal E_{\mathrm{start}}
 \le C\bigl(kh^r\sup_{0\le t\le 2k}\|u_t(t)\|_{H^r}
               +k^3\mathcal S_0\bigr).
\]
The constant is independent of $h$ and $k$.
\end{lemma}
\begin{proof}
Write $\theta^j=U_h^j-R_hu(t_j)$ and let $\tau_1,\tau_2$
be the starting residuals in the proof of
Lemma~\ref{lem:startup}; thus $\|\tau_j\|\le Ck^2\mathcal S_0$.
Subtracting the projected exact starting identities gives
\[
 \left(k^{-1}I+\tfrac12A_h\right)\theta^1
   =(P_h-R_h)D_1u(t_1)-P_h\tau_1,
\]
\[
 \left(\tfrac3{2k}I+A_h\right)\theta^2
   =\left(\tfrac2kI+2kb(k)A_h\right)\theta^1
       +(P_h-R_h)D_2u(t_2)-P_h\tau_2,
\]
where $\theta^0=0$.  By the $L^2$ projection's best approximation
property and \eqref{eq:Ritz},
\[
 \|(P_h-R_h)D_ju(t_j)\|
 \le Ch^r\sup_{0\le t\le2k}\|u_t(t)\|_{H^r},\qquad j=1,2.
\]
Positivity of $A_h$ gives
\[
 \left\|\left(k^{-1}I+\tfrac12A_h\right)^{-1}\right\|\le k,
\quad\left\|\left(\tfrac3{2k}I+A_h\right)^{-1}\right\|\le\tfrac{2k}{3},
\quad\left\|A_h\left(\tfrac3{2k}I+A_h\right)^{-1}\right\|\le1.
\]
It follows that
\[
 \|\theta^1\|\le C(kh^r\|u_t\|_{L^\infty H^r}+k^3\mathcal S_0),
\]
\[
 \|\theta^2\|
 \le(\tfrac43+2k|b(k)|)\|\theta^1\|
 +C(kh^r\|u_t\|_{L^\infty H^r}+k^3\mathcal S_0).
\]
Since $\mathcal E_{\mathrm{start}}=\|\theta^1\|+\|\theta^2\|$,
the claim follows.
\end{proof}

For piecewise linear finite elements on a quasi-uniform mesh and a sufficiently smooth solution, one may take $r=2$ in the $L^2$ norm.  If the weighted solution norms in Theorem~\ref{thm:FE} remain uniformly bounded and the starting values satisfy $\mathcal{E}_{\mathrm{start}}\le C(h^2+k^3)$ (as follows from Lemma~\ref{lem:FE-start} under its stated smoothness assumptions), then $\sup_{n\ge0}\norm{U_h^n-u(t_n)}\le C(h^2+k^3)$.

\section{Numerical experiments}\label{sec:numerics}

We here test the stability constant, temporal order, and spatial convergence.

\subsection{The stability constant}

For each fixed $\mu$, we generate $G_m(\mu)$ from \eqref{eq:Greenrec} and sum \eqref{eq:Lambda3mu} until the remaining tail is below machine precision.

On a coarse logarithmic grid in $\mu$ we first identify the neighborhood of the largest value; we then refine locally.  As a consistency check, the same procedure applied to BDF2 reproduces the value between $1.04$ and $1.05$ reported in \cite{ThomeeWahlbin1994}.

Table~\ref{tab:lambda3} shows representative BDF3 values.  The maximum found in the refined search is
\[
 \Lamthree\approx1.1501237996\qquad\text{at}\qquad \mu\approx7.42958394.
\]
The values close to one for small $\mu$ and for large $\mu$ agree with Proposition~\ref{prop:LambdaFinite}.  The computed maximum also lies inside the rigorous enclosure $23/20<\Lamthree<579/500$ from Theorem~\ref{thm:Lambda-enclosure}, and is only about $0.68\%$ below its analytic upper bound.

\begin{table}[htbp]
\caption{Representative values of $\Lambda_3(\mu)$.}\label{tab:lambda3}
\centering\footnotesize
\begin{tabular}{r r@{\qquad}r r}
\toprule
$\mu$ & $\Lambda_3(\mu)$ & $\mu$ & $\Lambda_3(\mu)$\\
\midrule
$0.01$ & $1.0000000$ & $1$ & $1.0186112$\\
$0.10$ & $1.0000000$ & $2$ & $1.1133169$\\
$0.50$ & $1.0000000$ & $5$ & $1.1294083$\\
$7.42958394$ & $1.1501238$ & $10$ & $1.1414634$\\
$20$ & $1.0938475$ & $100$ & $1.0273633$\\
\bottomrule
\end{tabular}
\end{table}

The numerical maximization is not used to justify Theorem~\ref{thm:stability}; the theorem relies on the exact supremum in Definition~\ref{def:Lambda3}.  It suggests the numerically suggested threshold $\beta\approx0.86947$, while Theorem~\ref{thm:Lambda-enclosure} supplies the rigorous explicit criterion $\beta<500/579\approx0.86356$.

\subsection{Convergence tests}

For the remaining experiments we take $\Omega=(0,1)$, $A=-\partial_{xx}$ with homogeneous Dirichlet boundary conditions, and
\[
 b(t)=\kappa e^{-at},\qquad u(x,t)=e^{-\rho t}\sin(\pi x).
\]
The forcing term is chosen from \eqref{eq:model}.  When $a\ne\rho$, it is
\[
 f(x,t)=\left[(\pi^2-\rho)e^{-\rho t}
 -\frac{\kappa\pi^2}{a-\rho}\bigl(e^{-\rho t}-e^{-at}\bigr)\right]\sin(\pi x).
\]
Unless otherwise stated, we use
\[
 \rho=0.7,\qquad a=2,\qquad \kappa=0.8.
\]
Then $\beta=\kappa/a=0.4$, and the condition \eqref{eq:kernel-sharp-condition} is safely satisfied.

\underline{Temporal convergence:} To isolate the temporal error, we evolve only the first eigenmode.  This is equivalent to solving the scalar problem with $\lambda=\pi^2$.  We report the maximum modal error over the whole interval $[0,1]$,
\[
 E_k=\max_{0\le n\le N}\abs{U^n-e^{-\rho t_n}}.
\]
Table~\ref{tab:timeconv} shows convergence to order three.  The pre-asymptotic rates are slightly below three because the error contains contributions from the two starting steps, but the rate approaches the predicted value monotonically as the grid is refined.

\begin{table}[htbp]
\caption{Temporal convergence for the manufactured solution on $[0,1]$.}\label{tab:timeconv}
\centering\footnotesize
\begin{tabular}{r c c c}
\toprule
$N$ & $k$ & $E_k$ & rate\\
\midrule
20   & $5.0000\times10^{-2}$ & $2.648911\times10^{-4}$ & --\\
40   & $2.5000\times10^{-2}$ & $4.186668\times10^{-5}$ & 2.662\\
80   & $1.2500\times10^{-2}$ & $5.975370\times10^{-6}$ & 2.809\\
160  & $6.2500\times10^{-3}$ & $8.016738\times10^{-7}$ & 2.898\\
320  & $3.1250\times10^{-3}$ & $1.048372\times10^{-7}$ & 2.935\\
640  & $1.5625\times10^{-3}$ & $1.352945\times10^{-8}$ & 2.954\\
1280 & $7.8125\times10^{-4}$ & $1.718601\times10^{-9}$ & 2.977\\
\bottomrule
\end{tabular}
\end{table}

\underline{Uniformity with respect to the final time:} We next fix $k=0.025$ and increase the final time while keeping all parameters unchanged.  We measure both the maximum error on $[0,T]$ and the error at the final time.  The results are shown in Table~\ref{tab:longtime}.  The maximum error is unchanged as $T$ grows from $1$ to $50$, while the final-time error decreases together with the exponentially decaying exact solution.  This behavior is exactly what the weighted estimate \eqref{eq:main-error} predicts.

\begin{table}[htbp]
\caption{Long-time behavior for fixed time-stepping $k=0.025$.}\label{tab:longtime}
\centering\footnotesize
\begin{tabular}{c c c}
\toprule
$T$ & $\max_{0\le t_n\le T}\abs{U^n-u(t_n)}$ & $\abs{U^N-u(T)}$\\
\midrule
1  & $4.186668\times10^{-5}$ & $5.102096\times10^{-6}$\\
5  & $4.186668\times10^{-5}$ & $6.429137\times10^{-7}$\\
10 & $4.186668\times10^{-5}$ & $2.161500\times10^{-8}$\\
20 & $4.186668\times10^{-5}$ & $2.001920\times10^{-11}$\\
50 & $4.186668\times10^{-5}$ & $1.518407\times10^{-20}$\\
\bottomrule
\end{tabular}
\end{table}

\underline{Spatial convergence:} Finally, we use continuous piecewise linear finite elements on a uniform partition of $(0,1)$.  We initialize with the Ritz projection of $u_0$, assemble the $L^2$ load vector for $P_hf$ exactly, and use the finite element counterparts of the two starting steps.  The time step is chosen as $k=h/2$, so the third-order time error is asymptotically smaller than the second-order $L^2$ spatial error.  The $L^2$ error at $T=1$ is listed in Table~\ref{tab:spaceconv}.  The observed rate is essentially two, in agreement with Theorem~\ref{thm:FE}.

\begin{table}[htbp]
\caption{Piecewise linear finite element convergence at $T=1$ with $k=h/2$.}\label{tab:spaceconv}
\centering\footnotesize
\begin{tabular}{r c c c}
\toprule
$1/h$ & $k$ & $L^2$ error & rate\\
\midrule
10  & $5.0000\times10^{-2}$ & $3.487757\times10^{-3}$ & --\\
20  & $2.5000\times10^{-2}$ & $8.711616\times10^{-4}$ & 2.001\\
40  & $1.2500\times10^{-2}$ & $2.175215\times10^{-4}$ & 2.002\\
80  & $6.2500\times10^{-3}$ & $5.433599\times10^{-5}$ & 2.001\\
160 & $3.1250\times10^{-3}$ & $1.357777\times10^{-5}$ & 2.001\\
\bottomrule
\end{tabular}
\end{table}

\appendix

\section{Proof of the quadrature estimate}
\label{app:quadrature}

We here prove Lemma~\ref{lem:quadrature}. The quadrature
\eqref{eq:quadrature-local} combines a Simpson starting
block on $[0,2k]$ with the classical three-step
Adams--Bashforth rule on $[t_\ell,t_{\ell+1}]$,
$\ell\ge2$; see \cite{HairerNorsettWanner1993}.
The error estimates follow from polynomial exactness
and integral Taylor remainders, as in the standard
analysis of interpolatory quadrature
\cite{DavisRabinowitz1984}.

Let $I_\ell g$ denote the quadratic interpolant
through $g(t_{\ell-2})$, $g(t_{\ell-1})$, and
$g(t_\ell)$. Direct integration gives
\[
 \int_{t_\ell}^{t_{\ell+1}}I_\ell g(s)\,ds
 =k\left(
 \frac{23}{12}g^\ell-\frac43g^{\ell-1}
 +\frac5{12}g^{\ell-2}
 \right).
\]
Set $J_\ell=[t_{\ell-2},t_{\ell+1}]$ and define
\[
 \mathcal E_\ell(g)
 :=\int_{t_\ell}^{t_{\ell+1}}
       \bigl(g(s)-I_\ell g(s)\bigr)\,ds.
\]
The functional $\mathcal E_\ell$ annihilates
polynomials of degree at most two. Moreover, since
the quadrature coefficients are fixed,
\[
 \norm{\mathcal E_\ell(v)}
 \le Ck\norm{v}_{C(J_\ell;H)}.
\]

For $g\in W^{3,1}(J_\ell;H)$, let
$T_{\ell,2}g$ be its quadratic Taylor polynomial
at $t_{\ell-2}$. The integral remainder formula
implies
\[
 \norm{g-T_{\ell,2}g}_{C(J_\ell;H)}
 \le Ck^2
 \int_{J_\ell}\norm{g^{(3)}(s)}\,ds.
\]
Since $\mathcal E_\ell(T_{\ell,2}g)=0$, we obtain
\begin{equation}\label{eq:local-peano3}
 \left\|
 \int_{t_\ell}^{t_{\ell+1}}
       \bigl(g(s)-I_\ell g(s)\bigr)\,ds
 \right\|
 \le Ck^3
 \int_{t_{\ell-2}}^{t_{\ell+1}}
       \norm{g^{(3)}(s)}\,ds.
\end{equation}
The Simpson starting block is also exact for
quadratic polynomials. Applying the same argument
on $[0,2k]$ yields
\begin{equation}\label{eq:start-peano3}
 \left\|
 \int_0^{2k}g(s)\,ds
 -\frac{k}{3}(g^0+4g^1+g^2)
 \right\|
 \le Ck^3
 \int_0^{2k}\norm{g^{(3)}(s)}\,ds.
\end{equation}
Each point of $(0,t_n)$ belongs to at most a
fixed number of the intervals $J_\ell$.
Summing the preceding estimates therefore gives
\[
 \norm{q^n(g)}
 \le Ck^3\int_0^{t_n}\norm{g^{(3)}(s)}\,ds,
\]
which proves \eqref{eq:quadrature-error} for $j=3$.

For $j=1$, we use exactness on constants.
Subtracting $g(t_{\ell-2})$ in the local error
functional and applying the fundamental theorem
of calculus gives
\[
 \norm{\mathcal E_\ell(g)}
 \le Ck
 \int_{t_{\ell-2}}^{t_{\ell+1}}
       \norm{g'(s)}\,ds.
\]
The Simpson starting block satisfies the analogous
estimate on $[0,2k]$. Summing with finite overlap
yields
\[
 \norm{q^n(g)}
 \le Ck\int_0^{t_n}\norm{g'(s)}\,ds.
\]
Thus \eqref{eq:quadrature-error} holds for both
$j=1$ and $j=3$, with constants independent of
$k$, $n$, and $g$.

It remains to verify the nonnegativity of the
assembled weights. For $n\ge5$, collecting the
contributions from the Simpson and Adams--Bashforth
blocks gives the first three coefficients $\frac34$, $\frac5{12}$, and $\frac43$. Each node with $3\le j\le n-3$ receives
$
 \frac{23}{12}-\frac43+\frac5{12}=1,
$
whereas the last two nodes receive $7/12$ and
$23/12$, respectively. These are precisely the
coefficients in \eqref{eq:quadrature-global}.
For $n=3,4$, direct collection gives
\[
 \sigma^3(g)
 =k\left(\frac34g^0+\frac94g^2\right)
\]
and
\[
 \sigma^4(g)
 =k\left(
 \frac34g^0+\frac5{12}g^1
 +\frac{11}{12}g^2+\frac{23}{12}g^3
 \right).
\]
Hence, all assembled quadrature weights are
nonnegative, and then we complete the proof.

\section{Bounds for the BDF3 resolvent constant}\label{app:Lambda-bound}

In this appendix, we prove the upper estimate in Theorem~\ref{thm:Lambda-enclosure}.

\subsection{A root parametrization and a positivity range}

We set
\[
 x=\frac{1}{11+6\mu},\qquad 0<x<\frac1{11}.
\]
Then $\eta=1-11x$, and the Green recurrence becomes
\begin{equation}\label{eq:x-Green}
 G_m=18xG_{m-1}-9xG_{m-2}+2xG_{m-3},\qquad m\ge3,
\end{equation}
with $G_0=1$, $G_1=18x$, and $G_2=324x^2-9x$.  In particular, every $G_m$ is a polynomial in $x$ with integer coefficients.

We first identify a spectral interval on which the Green sequence is nonnegative.

\begin{lemma}\label{lem:Green-positive}
If
\[
 \frac{256}{3915}\le x<\frac1{11},
\]
then $G_m(x)\ge0$ for every $m\ge0$.  Consequently, $\Lambda_3(\mu)=1$ throughout this interval.
\end{lemma}

\begin{proof}
Let $d=x^{-1}=11+6\mu$ and consider
\[
 d\zeta^3-18\zeta^2+9\zeta-2=0.
\]
Its discriminant equals
\[
   \triangle = -108(d^2-27d+189),
\]
and $d^2-27d+189>0$ for every real $d$.  Hence there is exactly one real root.  Since the polynomial is negative at $\zeta=0$ and positive at $\zeta=1$ for $d>11$, this root, denoted by $r$, belongs to $(0,1)$.  We write
\[
 D(r)=18r^2-9r+2,
\]
so the root equation gives
\begin{equation}\label{eq:x-r-map}
 d=\frac{D(r)}{r^3},\qquad x=\frac{r^3}{D(r)}.
\end{equation}
Moreover, we get
\[
 \frac{d}{dr}\frac{D(r)}{r^3}=-\frac{6(3r^2-3r+1)}{r^4}<0,
\]
so $r$ decreases strictly as $\mu$ increases.  At $r=8/15$, the formula \eqref{eq:x-r-map} gives $x=256/3915$.  Thus the range in the statement is equivalent to $r\ge8/15$.

Let the remaining two roots be $z=\rho e^{\mathrm{i}\theta}$ and $\bar z$.  The Vieta relations and \eqref{eq:x-r-map} give
\begin{align*}
   r \rho^2 = r z \bar z = 2x, \qquad r+ 2\rho \cos \theta  = r + z + \bar z = 18x,
\end{align*}
which further obtains
\begin{align}\label{rc}
 \rho^2=\frac{2r^2}{D(r)},\qquad
 \cos\theta=\frac{9r-2}{2\sqrt{2D(r)}}.
\end{align}
Since $r\geq 8/15 >1/2$ in the present range and $D(r)-2=9r(2r-1)$, we have $\rho\le r$.

The partial-fraction coefficients can be written explicitly in terms of
the three characteristic roots. Since the roots are $r$, $z=\rho e^{\mathrm{i}\theta}$, and
$\bar z=\rho e^{-\mathrm{i}\theta}$, the Green generating function is
\[
    \mathcal G(\xi)
    =
    \frac{1}{(1-r\xi)(1-z\xi)(1-\bar z\,\xi)}
    =
    \frac{A_0}{1-r\xi}
    +
    \frac{B_0}{1-z\xi}
    +
    \frac{\overline B_0}{1-\bar z\,\xi}.
\]
Expanding the three fractions into geometric series and comparing
coefficients of $\xi^m$ gives
\[
    G_m=A_0 r^m+B_0 z^m+\overline B_0\,\bar z^m.
\]
Multiplying the partial-fraction identity by $1-r\xi$ and setting
$\xi=r^{-1}$ yields
$A_0=r^2/[(r-z)(r-\bar z)]$.  Since
$z+\bar z=r(9r-2)/D(r)$ and
$z\bar z=\rho^2=2r^2/D(r)$, we obtain
\[
\begin{aligned}
    (r-z)(r-\bar z)
    &=
    r^2-r(z+\bar z)+z\bar z =
    \frac{r^2\bigl(D(r)-9r+4\bigr)}{D(r)}
    =
    \frac{6r^2E(r)}{D(r)},
\end{aligned}
\]
where $E(r):=3r^2-3r+1$.  Hence, we have $A_0=\frac{D(r)}{6E(r)}$. Similarly, multiplying by $1-z\xi$ and setting $\xi=z^{-1}$ gives
$B_0=z^2/[(z-r)(z-\bar z)]$, and therefore,
\[
    |B_0|^2
    =
    \frac{|z|^4}{|z-r|^2\,|z-\bar z|^2}.
\]
Here, $|z|^4=\rho^4=4r^4/D(r)^2$ and
$|z-r|^2=6r^2E(r)/D(r)$.  Moreover,
$z-\bar z=2\mathrm{i}\rho\sin\theta$, so that
\[
\begin{aligned}
    |z-\bar z|^2
    &=
    4\rho^2(1-\cos^2\theta) =
    \frac{r^2\bigl(8D(r)-(9r-2)^2\bigr)}{D(r)^2} =
    \frac{3r^2(21r^2-12r+4)}{D(r)^2},
\end{aligned}
\]
where we used \eqref{rc}. Consequently,
\begin{equation}\label{eq:AB-formulas}
    |B_0|^2=
    \frac{2D(r)}
    {9E(r)(21r^2-12r+4)}.
\end{equation}
The denominators in $A_0$ and \eqref{eq:AB-formulas} are strictly positive. Furthermore, we have
\[
 A_0^2-4\abs{B_0}^2=
 \frac{D(r)Q(r)}{12E(r)^2(21r^2-12r+4)},
\]
where
\[
 Q(r)=126r^4-135r^3+42r^2+12r-8.
\]
At $r=8/15$, we have $Q(8/15)=344/5625>0$.  To see that $Q(r)$ remains positive, we write
\[
 Q'(r)=3S(r),\qquad S(r)=168r^3-135r^2+28r+4.
\]
For $r\ge8/15$, we yield
\[
 S(8/15)=\frac{6772}{1125}>0,
 \qquad S'(8/15)=\frac{684}{25}>0,
 \qquad S''(r)=1008r-270>0.
\]
Thus $S'$ and $S$ are increasing and positive on this interval, and hence $Q(r)>0$.  We conclude that $A_0 \ge2\abs{B_0}$.  Combining this with $\rho\le r$ yields
\[
 G_m\ge A_0 r^m-2\abs{B_0}\rho^m\ge (A_0-2\abs{B_0})r^m\ge0,
\]
for every $m\ge0$.  Finally, the signed-mass identity \eqref{eq:signedmass} gives $\eta\sum_{m\ge0}G_m=1$, so Definition~\ref{def:Lambda3} implies $\Lambda_3(\mu)=1$.
\end{proof}

\subsection{Exact Bernstein enclosures}

We next introduce a polynomial enclosure that will be used to control
the finite part of the Green series uniformly on rational subintervals.
Its main advantage is that, when the polynomial coefficients and the
interval endpoints are rational, the resulting bounds can be evaluated
entirely by exact rational arithmetic.

\begin{lemma}\label{lem:Bernstein}
Let $p(x)=\sum_{q=0}^{\nu}c_qx^q$ and let $I=[a,b]$.  Under the affine
change of variables $x=a+(b-a)t$, $0\le t\le1$, write
\[
 p(a+(b-a)t)
 =
 \sum_{j=0}^{\nu}\mathcal{B}_j(p;I)
 \binom{\nu}{j}t^j(1-t)^{\nu-j}.
\]
Then, for every $x\in I$,
\[
 \min_{0\le j\le\nu}\mathcal{B}_j(p;I)
 \le p(x)\le
 \max_{0\le j\le\nu}\mathcal{B}_j(p;I).
\]
If $a$, $b$, and all coefficients $c_q$ are rational, then every
Bernstein coefficient is rational and is given explicitly by
\begin{equation}\label{eq:Bernstein-coeff}
 \mathcal{B}_j(p;I)=
 \sum_{\ell=0}^{j}
 \frac{\binom{j}{\ell}}{\binom{\nu}{\ell}}
 (b-a)^\ell
 \sum_{q=\ell}^{\nu}
 c_q\binom{q}{\ell}a^{q-\ell}.
\end{equation}
\end{lemma}

\begin{proof}
We set the Bernstein basis functions
\[
 B_{j,\nu}(t)
 :=
 \binom{\nu}{j}t^j(1-t)^{\nu-j},
 \qquad 0\le j\le\nu.
\]
For $0\le t\le1$, one has $B_{j,\nu}(t)\ge0$ and
$\sum_{j=0}^{\nu}B_{j,\nu}(t)=1$.  Hence
$p(a+(b-a)t)$ is a convex combination of the coefficients
$\mathcal{B}_j(p;I)$, which immediately gives the stated enclosure.

To obtain \eqref{eq:Bernstein-coeff}, first expand
$p(a+(b-a)t)$ in powers of $t$.  Using the binomial formula, we have
\[
 p(a+(b-a)t)
 =
 \sum_{\ell=0}^{\nu}
 \left[
 (b-a)^\ell
 \sum_{q=\ell}^{\nu}
 c_q\binom{q}{\ell}a^{q-\ell}
 \right]t^\ell,
\]
where $t^\ell$ can be rewritten in the Bernstein basis as $t^\ell =\sum_{j=\ell}^{\nu} \binom{j}{\ell}/\binom{\nu}{\ell} \,  B_{j,\nu}(t)$.
Substitution and comparison of Bernstein coefficients give
\eqref{eq:Bernstein-coeff}.  The rationality assertion is immediate
from the explicit formula.
\end{proof}

For an interval $I$, define
\[
 \mathfrak B_I[p]
 :=
 \max_{0\le j\le\nu}
 \abs{\mathcal{B}_j(p;I)}.
\]
Lemma~\ref{lem:Bernstein} then yields the uniform bound
$\abs{p(x)}\le\mathfrak B_I[p]$ for every $x\in I$.  Thus a polynomial
can be bounded on the whole interval by finitely many explicitly
computable Bernstein coefficients.  In particular, on rational
intervals these bounds are exact rational quantities rather than
floating-point estimates.

We apply this construction to the Green coefficients.  Recall that
$x=(11+6\mu)^{-1}$, so
$\eta(\mu)=\frac{6\mu}{11+6\mu}=1-11x$.
Moreover, the recurrence
\[
 G_m(x)=18xG_{m-1}(x)-9xG_{m-2}(x)+2xG_{m-3}(x)
\]
with $G_0=1$, $G_1=18x$, and $G_2=324x^2-9x$ shows that every
$G_m(x)$ is a polynomial in $x$ with integer coefficients.  We therefore
set
\[
 P_m(x):=(1-11x)G_m(x).
\]
Since $0<x<1/11$, the factor $1-11x$ is positive, and hence
\[
 \Lambda_3(\mu)
 =
 (1-11x)\sum_{m=0}^{\infty}\abs{G_m(x)}
 =
 \sum_{m=0}^{\infty}\abs{P_m(x)}.
\]
Consequently, for every rational interval $I\subset(0,1/11)$ and every
finite $N$, we get
\[
 \sum_{m=0}^{N}\abs{P_m(x)}
 \le
 \sum_{m=0}^{N}\mathfrak B_I[P_m],
 \qquad x\in I.
\]
The right-hand side is an exact rational upper bound for the finite part
of the resolvent series.  In the next two subsections this finite
Bernstein enclosure will be combined with separate analytic estimates
for the tail $\sum_{m>N}\abs{P_m(x)}$.

\subsection{The intermediate range}

We next consider
\[
 \frac1{32}\le x\le\frac{256}{3915}.
\]
Introduce the companion matrix
\[
 \widetilde{C}(x)=
 \begin{pmatrix}
 18x&-9x&2x\\
 1&0&0\\
 0&1&0
 \end{pmatrix},
 \qquad
 y_m=(G_m,G_{m-1},G_{m-2})^T.
\]
Then $y_{m+1}=\widetilde{C}(x)y_m$.  Applying Lemma~\ref{lem:Bernstein} entrywise to the polynomial matrices $\widetilde{C}(x)^s$ over the entire intermediate interval gives the following exact rational bounds:
\begin{center}
\begin{tabular}{c|cccccc}
$s$&0&1&2&3&4&5\\ \hline
upper bound for $\Big\|\widetilde{C}(x)^s\Big\|_1$
&$1$&$109/50$&$149/50$&$47/20$&$131/100$&$57/100$
\end{tabular}
\end{center}
In particular, $\|\widetilde{C}^5\|_1<57/100$ and $\sum_{s=0}^4\|\widetilde{C}^s\|_1<10$.  Therefore, for every $N\ge2$ and $m=N+5j+s$ with $j\geq 0$ and $s=0,1,2,3,4$, we have
\begin{align*}
    II:= \sum_{m=N}^{\infty}\abs{G_m} \leq \sum_{m=N}^{\infty}\|y_m\|_1 = \sum_{j=0}^{\infty}\sum_{s=0}^4 \|y_{N+5j+s}\|_1,
\end{align*}
which combines $y_{m+1}=\widetilde{C} \, y_m$ and $\sum_{j=0}^{\infty}(57/100)^j = \frac{100}{43}$ to give
\begin{equation}\label{eq:intermediate-tail}
 II
 \le\sum_{j=0}^{\infty}\sum_{s=0}^4
 \|\widetilde{C}^s (\widetilde{C}^{5j}) y_N\|_1 \leq \|y_N\|_1 \sum_{j=0}^{\infty} \|(\widetilde{C}^5)^j\|_1 \sum_{s=0}^4 \|\widetilde{C}^s\|_1
 <\frac{1000}{43}\|y_N\|_1.
\end{equation}

We partition the interval into $48$ equal rational subintervals
\[
 I_j^{\rm M}=\left[\frac1{32}+\frac{j}{48}\left(\frac{256}{3915}-\frac1{32}\right), \,
 \frac1{32}+\frac{j+1}{48}\left(\frac{256}{3915}-\frac1{32}\right)\right],
 \qquad 0\le j\le47.
\]
We set $B_{m,j}^{\rm M}=\mathfrak B_{I_j^{\rm M}}[P_m]$.  Combining \eqref{eq:intermediate-tail} with $N=25$ gives
\[
 \Lambda_3(\mu)
 \le \sum_{m=0}^{24}B_{m,j}^{\rm M}
 +\frac{1000}{43}\left(B_{23,j}^{\rm M}+B_{24,j}^{\rm M}+B_{25,j}^{\rm M}\right),
 \qquad x\in I_j^{\rm M}.
\]
Every term on the right is an exact rational number computable from \eqref{eq:Bernstein-coeff}.  Evaluating the $48$ finite rational expressions by integer arithmetic gives the slightly sharper certificate
\begin{equation}\label{eq:intermediate-certificate}
 \max_{0\le j\le47}
 \left\{\sum_{m=0}^{24}B_{m,j}^{\rm M}
 +\frac{1000}{43}\left(B_{23,j}^{\rm M}+B_{24,j}^{\rm M}+B_{25,j}^{\rm M}\right)\right\}
 <\frac{2863}{2500}<\frac{573}{500}.
\end{equation}
The largest enclosure occurs for $j=6$, whose interval is
$I_6^{\rm M}= \left[\frac{35597}{1002240}, \, \frac{217859}{6013440}\right]$.
The location is reported only to make the finite verification easier to reproduce; inequality \eqref{eq:intermediate-certificate} itself is an exact rational comparison.

\subsection{The high-frequency range}

It remains to consider $0<x\le1/32$.  Here a direct scalar tail estimate is sharper.  Let
\[
 T_N=\sum_{m=N}^{\infty}\abs{G_m}.
\]
From \eqref{eq:x-Green}, we get
\[
 T_N\le x\bigl(18T_{N-1}+9T_{N-2}+2T_{N-3}\bigr).
\]
Since
\[
 T_{N-1}=\abs{G_{N-1}}+T_N,
\]
and similarly for $T_{N-2}$ and $T_{N-3}$, we obtain
\begin{equation}\label{eq:large-tail}
 (1-29x)T_N
 \le x\left(29\abs{G_{N-1}}+11\abs{G_{N-2}}+2\abs{G_{N-3}}\right).
\end{equation}
The coefficient on the left is positive because $x\le1/32$.

Partition this interval into
\[
 I_j^{\rm L}=\left[\frac{j}{2048},\frac{j+1}{2048}\right],
 \qquad 0\le j\le63,
\]
and put $B_{m,j}^{\rm L}=\mathfrak B_{I_j^{\rm L}}[P_m]$ and $b_j=(j+1)/2048$.  Multiplying \eqref{eq:large-tail} by $1-11x$ and taking $N=30$ yields
\[
 \Lambda_3(\mu)\le
 \sum_{m=0}^{29}B_{m,j}^{\rm L}
 +\frac{b_j}{1-29b_j}
 \left(29B_{29,j}^{\rm L}+11B_{28,j}^{\rm L}+2B_{27,j}^{\rm L}\right),
 \qquad x\in I_j^{\rm L}.
\]
Again all quantities are rational.  Applying \eqref{eq:Bernstein-coeff} gives
\begin{equation}\label{eq:large-certificate}
 \max_{0\le j\le63}
 \left[\sum_{m=0}^{29}B_{m,j}^{\rm L}
 +\frac{b_j}{1-29b_j}
 \left(29B_{29,j}^{\rm L}+11B_{28,j}^{\rm L}+2B_{27,j}^{\rm L}\right)\right]
 <\frac{2893}{2500}<\frac{579}{500}.
\end{equation}
The largest enclosure is attained at $j=36$, namely $I_{36}^{\rm L}=[9/512,37/2048]$.

Combining Lemma~\ref{lem:Green-positive}, \eqref{eq:intermediate-certificate}, and \eqref{eq:large-certificate}, we obtain
\[
 \Lambda_3<\max\left\{1,\frac{573}{500},\frac{579}{500}\right\}=\frac{579}{500},
\]
which completes the proof of the upper bound in Theorem~\ref{thm:Lambda-enclosure}.

\begin{remark}
The finite verifications in \eqref{eq:intermediate-certificate} and \eqref{eq:large-certificate} are exact rational comparisons.  Starting from the integer-coefficient recurrence \eqref{eq:x-Green}, one forms $P_m=(1-11x)G_m$, converts each $P_m$ to Bernstein form by \eqref{eq:Bernstein-coeff}, takes the largest absolute Bernstein coefficient on each rational subinterval, and finally performs the displayed rational additions and cross-multiplications.  For orientation only, the largest two rational enclosures have decimal values approximately $1.145170316$ and $1.157152284$, respectively; these decimal representations play no role in the proof.
\end{remark}

\section{Weighted resolvent estimates}\label{app:weighted}

We prove Lemma~\ref{lem:weighted-resolvent}.  The difficulty is uniformity in the coupled parameters $k$, $\lambda$, and $\mu=k\lambda$.  One root approaches 1 when $\mu\to0$, while for large $\mu$ all recurrence coefficients are small.  We therefore divide the positive spectrum into three fixed regions and prove a separate lemma on each of them.

Throughout this appendix let
\[
 \varpi=e^{\gamma k}\ge1.
\]
The weighted Green sum is $\eta(\mu)\sum_{m\ge0}\varpi^m\abs{G_m(\mu)}$. The BDF3 symbol $\delta(z)=\frac{11}{6}-3z+\frac32z^2-\frac13z^3$ satisfies
\[
 \left(\frac{\delta(z)}k+\lambda\right)^{-1}
 =\frac{6k}{11+6\mu}\,\mathcal G_\mu(z),\qquad \mu=k\lambda.
\]

\subsection{The small spectral region}

\begin{lemma}\label{lem:weighted-small}
Let $\lambda_0>0$ and $\varepsilon>0$.  There exist $\mu_s\in(0,1/2]$, $\gamma_s>0$, $k_s>0$, and $C_s>0$ such that, whenever
\[
 0<\mu=k\lambda\le\mu_s,\qquad \lambda\ge\lambda_0,
 \qquad 0\le\gamma\le\gamma_s,
 \qquad 0<k\le k_s,
\]
one has
\begin{align}\label{etg}
  \eta(\mu)\sum_{m=0}^{\infty}\varpi^m\abs{G_m(\mu)}\le 1+\varepsilon
\end{align}
and
\begin{align}\label{etg2}
 \sup_{m\ge0}\varpi^m\abs{G_m(\mu)}\le C_s.
\end{align}
Every homogeneous solution satisfies
\[
 \abs{h^n}\le C_s e^{-\gamma t_n}
 \bigl(\abs{h^0}+\abs{h^1}+\abs{h^2}\bigr).
\]
\end{lemma}

\begin{proof}
Let
\[
 F(r,\mu)=(11+6\mu)r^3-18r^2+9r-2.
\]
At $(r,\mu)=(1,0)$, $F=0$, $F_r=6$, and $F_\mu=6$.  The implicit-function theorem therefore gives a real analytic principal root $r_1(\mu)$ with
\[
 r_1(0)=1,\qquad r_1'(0)=-1.
\]
After reducing $\mu_s$, we have
\begin{align*}
    r_1(\mu) = 1 - \mu + O(\mu^2) \leq  1 - \mu + C\mu^2 \leq 1 - \mu + \frac{1}{4} \mu = 1-\frac{3}{4}\mu,
\end{align*}
which follows the basic inequality $1-x \leq e^{-x}$ with $x\geq 0$, then
\begin{equation}\label{eq:small-root-decay}
 0<r_1(\mu)\le e^{-3\mu/4},\qquad 0\le\mu\le\mu_s.
\end{equation}
At $\mu=0$, the two parasitic roots are $(7\pm \mathrm{i}\sqrt{39})/22$ and have modulus $\sqrt{2/11}$.  By continuity there is a fixed $\rho_s<1$ such that
\[
 \max\{\abs{r_2(\mu)},\abs{r_3(\mu)}\}\le\rho_s,
 \qquad 0\le\mu\le\mu_s.
\]
For sufficiently small $\mu_s$ the three roots remain distinct.  Hence, the partial-fraction coefficients are continuous and uniformly bounded, and we use \eqref{eq:Green-root-decomp} to write
\begin{equation}\label{GGA}
  G_m(\mu)=A_1(\mu)r_1(\mu)^m+R_m(\mu),
 \qquad \abs{R_m(\mu)}\le C\rho_s^m
\end{equation}
with $A_1(0)=11/6$.  Also, for $\mu<1$, we have
\begin{equation}\label{eq:eta-small}
 \eta(\mu)=\frac{6}{11}\mu+O(\mu^2).
\end{equation}
Since $\mu=k\lambda$ and $\lambda\ge\lambda_0$, we have $k\le\mu/\lambda_0$.  First choose $\gamma_s\le\lambda_0/8$.  From \eqref{eq:small-root-decay},
\begin{equation}\label{qr1}
 \varpi r_1(\mu)\le
 \exp\left(\frac{\gamma}{\lambda_0}\mu-\frac34\mu\right)
 \le e^{-5\mu/8}.
\end{equation}
Hence, after reducing $\mu_s$ if necessary,
\[
 1-qr_1(\mu)\ge1-e^{-5\mu/8}\ge\frac{5}{16}\mu.
\]
This already gives uniform summability of the principal mode after multiplication by $\eta(\mu)$.  To obtain the sharp limit needed below, write $\vartheta=\gamma/\lambda$.  Then $0\le\vartheta\le\gamma_s/\lambda_0$ and $\varpi=e^{\vartheta\mu}=1+\vartheta\mu + O(\mu^2)$.  The expansions $r_1(\mu)=1-\mu+O(\mu^2)$, $A_1(\mu)=11/6+O(\mu)$, and \eqref{eq:eta-small}, uniformly for bounded $\vartheta$, give
\[
   \eta(\mu)\sum_{m=0}^{\infty}\varpi^m\abs{A_1(\mu)} r_1(\mu)^m = \frac{\eta(\mu)\abs{A_1(\mu)}}{1-e^{\vartheta\mu}r_1(\mu)}
 =\frac{1+O(\mu)}{1-\vartheta+O(\mu)} := \frac{1+\delta_1(\mu,\vartheta)}{1-\vartheta+\delta_2(\mu,\vartheta)}.
\]
After fixing the upper bound $\gamma_s\le\lambda_0/8$, decrease $\gamma_s$ further so that
\begin{align*}
     \frac{1}{1-\vartheta} \leq \frac{1}{1-\gamma_s/\lambda_0}< 1+ \frac{\varepsilon}{4},
\end{align*}
and then choose $\mu_s$ so small that
\begin{equation*}
   \Big| \frac{1+\delta_1(\mu,\vartheta)}{1-\vartheta+\delta_2(\mu,\vartheta)} - \frac{1}{1-\vartheta} \Big| \leq \frac{\varepsilon}{4}.
\end{equation*}
Thus, for $0<\mu\leq \mu_s$ and $0\leq \vartheta \leq \gamma_s/\lambda_0$, we have
\begin{align}\label{eq:small-principal-bound}
   \eta(\mu)\sum_{m=0}^{\infty}\varpi^m\abs{A_1(\mu)} r_1(\mu)^m \le \frac{1}{1-\vartheta} + \frac{\varepsilon}{4} \leq 1+ \frac{\varepsilon}{2}.
\end{align}
The parasitic part satisfies, after taking $k_s$ small enough that $\varpi\rho_s<\widetilde\rho_s<1$,
\[
 \eta(\mu)\sum_{m=0}^{\infty}\varpi^m\abs{R_m(\mu)}
 \le \frac{C\eta(\mu)}{1-\widetilde\rho_s} \leq \frac{C\mu}{1-\widetilde\rho_s} \leq  \frac{C\mu_s}{1-\widetilde\rho_s} \leq \frac{\varepsilon}{2},
\]
after decreasing $\mu_s$ further; together with \eqref{eq:small-principal-bound}, proves the weighted $\ell^1$ estimate \eqref{etg}.

We next prove the pointwise estimate \eqref{etg2}.  Since the
partial-fraction coefficient $A_1(\mu)$ depends continuously on $\mu$
and $[0,\mu_s]$ is compact, there exists a constant $M_s>0$ such that
\[
    \sup_{0\le\mu\le\mu_s}|A_1(\mu)|\le M_s.
\]
Using \eqref{GGA}, we obtain
\[
\begin{aligned}
    \varpi^m|G_m(\mu)|
    &\le
    |A_1(\mu)|(\varpi r_1(\mu))^m
    +C(\varpi\rho_s)^m.
\end{aligned}
\]
By the choices made above, then
\[
    \varpi r_1(\mu)\le e^{-5\mu/8}\le1,
    \qquad
    \varpi\rho_s\le\widetilde\rho_s<1.
\]
Consequently,
\[
\begin{aligned}
    \varpi^m|G_m(\mu)|
    \le M_s(\varpi r_1(\mu))^m
    +C\widetilde\rho_s^m \le M_s+C,
\end{aligned}
\]
which proves \eqref{etg2}.

It remains to establish the exponentially weighted bound for the
homogeneous recurrence.  Let $h^n$ satisfy
\[
    h^n-a_1(\mu)h^{n-1}
    +a_2(\mu)h^{n-2}
    -a_3(\mu)h^{n-3}=0,
    \qquad n\ge3.
\]
Since the three characteristic roots remain distinct on
$[0,\mu_s]$, the solution admits the representation
\begin{equation}\label{eq:small-homogeneous-root-rep}
    h^n
    =
    c_1(\mu)r_1(\mu)^n
    +c_2(\mu)r_2(\mu)^n
    +c_3(\mu)r_3(\mu)^n,
    \qquad n\ge0,
\end{equation}
where the coefficients $c_j(\mu)$ are determined by the starting
values $h^0,h^1,h^2$.  More precisely,
\[
    V(\mu)
    \begin{pmatrix}
        c_1(\mu)\\
        c_2(\mu)\\
        c_3(\mu)
    \end{pmatrix}
    =
    \begin{pmatrix}
        h^0\\
        h^1\\
        h^2
    \end{pmatrix},
    \qquad
    V(\mu)
    :=
    \begin{pmatrix}
        1&1&1\\
        r_1(\mu)&r_2(\mu)&r_3(\mu)\\
        r_1(\mu)^2&r_2(\mu)^2&r_3(\mu)^2
    \end{pmatrix}.
\]
The determinant of this Vandermonde matrix is
\[
    \det V(\mu)
    =
    \bigl(r_2(\mu)-r_1(\mu)\bigr)
    \bigl(r_3(\mu)-r_1(\mu)\bigr)
    \bigl(r_3(\mu)-r_2(\mu)\bigr).
\]
The roots are distinct for every $\mu\in[0,\mu_s]$.  By continuity
and compactness, after reducing $\mu_s$ if necessary, there exists
$\delta_s>0$ such that
\[
    \min_{0\le\mu\le\mu_s}
    \min_{i\ne j}
    |r_i(\mu)-r_j(\mu)|
    \ge\delta_s.
\]
Hence, $\det V(\mu)$ stays uniformly away from zero.  Since the entries
of $V(\mu)$ are also uniformly bounded on $[0,\mu_s]$, its inverse is
uniformly bounded: there exists $C_V>0$ such that
\[
    \sup_{0\le\mu\le\mu_s}
    \|V(\mu)^{-1}\|_1
    \le C_V.
\]
It follows that
\begin{equation*}
    |c_1(\mu)|+|c_2(\mu)|+|c_3(\mu)|
    \le
    C_V\bigl(|h^0|+|h^1|+|h^2|\bigr).
\end{equation*}
Multiplying \eqref{eq:small-homogeneous-root-rep} by $\varpi^n$ and using
the triangle inequality gives
\[
\begin{aligned}
    \varpi^n|h^n|
    &\le
    |c_1(\mu)|(\varpi r_1(\mu))^n
    +|c_2(\mu)|(\varpi|r_2(\mu)|)^n
    +|c_3(\mu)|(\varpi|r_3(\mu)|)^n.
\end{aligned}
\]
The principal-root estimate \eqref{qr1} gives $\varpi r_1(\mu)\le e^{-5\mu/8}<1$ for
$\mu>0$, while
\[
    \varpi|r_2(\mu)|,\ \varpi|r_3(\mu)|
    \le \varpi\rho_s
    \le\widetilde\rho_s<1.
\]
Therefore, we have
\[
\begin{aligned}
    \varpi^n|h^n|
    &\le
    |c_1(\mu)|+|c_2(\mu)|+|c_3(\mu)| \le
    C_V\bigl(|h^0|+|h^1|+|h^2|\bigr).
\end{aligned}
\]
Since $\varpi=e^{\gamma k}$ and $t_n=nk$, we have
$\varpi^n=e^{\gamma t_n}$. Thus,
\[
    |h^n|
    \le
    C_Ve^{-\gamma t_n}
    \bigl(|h^0|+|h^1|+|h^2|\bigr).
\]
Taking $C_s=\max\{M_s+C,C_V\}$ completes the two uniform estimates.
\end{proof}

\subsection{The compact intermediate spectral region}

\begin{lemma}\label{lem:weighted-middle}
Fix $0<\mu_s<M<\infty$ and $\varepsilon>0$.  There exist $\gamma_m>0$, $k_m>0$, and $C_m>0$ such that, if
\[
 \mu_s\le\mu\le M,\qquad 0\le\gamma\le\gamma_m,
 \qquad 0<k\le k_m,
\]
then
\[
 \eta(\mu)\sum_{m=0}^{\infty}\varpi^m\abs{G_m(\mu)}
 \le\Lamthree+\varepsilon, \quad \sup_{m\geq 0} \varpi^m\abs{G_m(\mu)}\le C_m,
\]
and the corresponding homogeneous solutions satisfy the exponentially weighted estimate in Lemma~\ref{lem:weighted-resolvent}.
\end{lemma}

\begin{proof}
By Lemma~\ref{lem:schur} and compactness, there is $\rho_m<1$ such that every characteristic root satisfies $\abs{r_j(\mu)}\le\rho_m$ for $\mu\in[\mu_s,M]$.  The discriminant of $(11+6\mu)r^3-18r^2+9r-2$ is $-108(36\mu^2-30\mu+13)$, which is strictly negative for every $\mu>0$.  Thus the three roots are simple and their partial-fraction coefficients are uniformly bounded on the compact interval.  Hence,
\begin{equation}\label{eq:middle-G-bound}
 \abs{G_m(\mu)}\le C\rho_m^m,
 \qquad \mu_s\le\mu\le M,
\end{equation}
with $C$ independent of $m$ and $\mu$.

Choose $\varpi_*>1$ so close to one that $\varpi_*\rho_m<1$.  If $\varpi\le \varpi_*$, then by the mean-value theorem,
\[
 \varpi^m-1\le m \varpi_*^{m-1}(\varpi-1).
\]
Using \eqref{eq:middle-G-bound} and $0<\eta(\mu)<1$ gives
\begin{align*}
 &\eta(\mu)\sum_{m=0}^{\infty}(\varpi^m-1)\abs{G_m(\mu)}\le C(\varpi-1)\sum_{m=1}^{\infty}m \varpi_*^{m-1}\rho_m^m
 \le C_*(\varpi-1).
\end{align*}
Since $\varpi-1=e^{\gamma k}-1$, we may choose $\gamma_m$ and $k_m$ so that $\varpi\le \varpi_*$ and $C_*(\varpi-1)\le\varepsilon$.  Therefore, we get
\[
 \eta(\mu)\sum_{m=0}^{\infty}\varpi^m\abs{G_m(\mu)}
 \le\LamthreeMu+\varepsilon\le\Lamthree+\varepsilon.
\]
The pointwise estimate follows directly from \eqref{eq:middle-G-bound}.  The same uniform root representation applies to the homogeneous solution and, because $\varpi\rho_m<1$, gives the required decay after multiplication by $e^{-\gamma t_n}$.
\end{proof}

\subsection{The large spectral region}

\begin{lemma}\label{lem:weighted-large}
Let $M=40$.  For every $\varepsilon>0$ there exist $\gamma_\ell>0$, $k_\ell>0$, and $C_\ell>0$ such that, if
\[
 \mu\ge M,
 \qquad 0\le\gamma\le\gamma_\ell,
 \qquad 0<k\le k_\ell,
\]
then it holds that
\[
 \eta(\mu)\sum_{m=0}^{\infty}\varpi^m\abs{G_m(\mu)}
 \le\Lamthree+\varepsilon, \quad \sup_{m\geq 0} \varpi^m\abs{G_m(\mu)}\le C_\ell,
\]
and every homogeneous solution has the weighted decay stated in Lemma~\ref{lem:weighted-resolvent}.
\end{lemma}

\begin{proof}
For $\mu\ge40$, the recurrence coefficients satisfy
\[
 a_1+a_2+a_3=\frac{29}{11+6\mu}\le\frac{29}{251}.
\]
The weighted sequence $H_m=\varpi^mG_m$ has generating function
\[
 \sum_{m=0}^{\infty}H_m z^m=
 \frac{1}{1-a_1qz+a_2(qz)^2-a_3(qz)^3}.
\]
Let
\[
 S_{\varpi}(\mu)=a_1\varpi+a_2\varpi^2+a_3\varpi^3.
\]
Choose $\gamma_\ell k_\ell$ so small that, uniformly for $\mu\ge40$, then $S_{\varpi}(\mu) \leq \frac{29\varpi^3 }{11+6\mu} \leq \frac{3}{25}$. Expanding the generating function as a Neumann series in the Wiener algebra gives
\begin{equation}\label{eq:large-Wiener}
 \sum_{m=0}^{\infty}\varpi^m\abs{G_m(\mu)}
 \le\frac{1}{1-S_{\varpi}(\mu)}\le\frac{25}{22}.
\end{equation}
Since $\eta<1$ and $25/22<23/20<\Lamthree$ by Theorem~\ref{thm:Lambda-enclosure}, \eqref{eq:large-Wiener} is already stronger than the required $\Lamthree+\varepsilon$ estimate.  It also gives the pointwise Green bound.

For a homogeneous solution define $H^n=\varpi^n\abs{h^n}$.  The recurrence and the triangle inequality imply that
\[
 H^n\le a_1\varpi H^{n-1}+a_2\varpi^2H^{n-2}+a_3\varpi^3H^{n-3}
 \le S_{\varpi}\max_{n-3\le j\le n-1}H^j.
\]
Because $S_{\varpi}<1$, the induction gives $H^n\le C(H^0+H^1+H^2)$, which is equivalent to the asserted exponential decay.
\end{proof}

\subsection{Proof of Lemma~\ref{lem:weighted-resolvent}}

We now assemble the three regimes.  Apply Lemma~\ref{lem:weighted-small} with $\varepsilon$ replaced by $\varepsilon/2$ and let $\mu_s$ be the resulting small cutoff.  Use Lemma~\ref{lem:weighted-middle} on the compact interval $[\mu_s,40]$, and Lemma~\ref{lem:weighted-large} for $\mu\ge40$.  Choose
\[
 \gamma_0=\min\{\gamma_s,\gamma_m,\gamma_\ell\},
 \qquad
 k_0=\min\{k_s,k_m,k_\ell\},
\]
and let $C$ be the maximum of the three regional constants.  Since $\Lamthree>1$, the small-region estimate $1+\varepsilon/2$ is bounded by $\Lamthree+\varepsilon$.  The three regional lemmas therefore give \eqref{eq:weighted-Lambda}, \eqref{eq:weighted-pointwise}, and \eqref{eq:homogeneous-bound} uniformly for every $\lambda\ge\lambda_0$.  This proves Lemma~\ref{lem:weighted-resolvent}.

\end{document}